\documentclass[reqno]{amsproc}
\usepackage{tikz}

\usepackage{amstext,amsmath,amssymb,amsfonts}
\usepackage[latin1]{inputenc}
\usepackage{epsfig}
\usepackage{hyperref}
\usepackage{color}

\usepackage{latexsym}

\newtheorem{lemma}{Lemma}
\newtheorem{corollary}{Corollary}

\newtheorem{definition}{Definition}
\newtheorem{theorem}{Theorem}

\newtheorem{proposition}{Proposition}

\newcommand{\vspi}{\vspace{0.4cm}}

\newcommand{\bea}{\begin{eqnarray}}
\newcommand{\eea}{\end{eqnarray}}
\newcommand{\beq}{\begin{equation}}
\newcommand{\eeq}{\end{equation}}
\newcommand{\enn}{\nonumber \end{equation}}

\newcommand{\tw}{ {\rm t}}
\newcommand{\unt}{ {\rm unt}}
 \newcommand{\rk}{{\rm r}\,}

 \newcommand{\cG}{\mathcal{G}}
 \newcommand{\cV}{\mathcal{V}}
 \newcommand{\cE}{\mathcal{E}}
 \newcommand{\cF}{\mathcal{F}}

 \newcommand{\cR}{\mathcal{R}}
  \newcommand{\cQ}{\mathcal{Q}}

 \newcommand{\cC}{\mathcal{C}}

\newcommand{\bG}{\partial{\mathcal{G}}}
\newcommand{\bV}{{\mathcal{V}}_{\partial}}
\newcommand{\bE}{{\mathcal{E}}_{\partial}}
\newcommand{\bC}{{\mathcal{C}}_{\partial}}

 \newcommand{\mf}{\mathfrak{f}}

 \newcommand{\sset}{\Subset}

 \newcommand{\inter}{{\rm int}}
 \newcommand{\ext}{{\rm ext}}

\title[Universality for polynomials invariants for ribbon graph with HRs]{ {Universality for polynomial invariants} \\{for ribbon 
graphs with half-ribbons}\footnote{Preprint: ICMPA-MPA/2013/10}}

\author{R\'emi C. Avohou}
\address[R.C.A.]{Humboldt-Universit{\"a}t zu Berlin, Institut f{\"u}r Mathematik und Institut f{\"u}r Physik, Rudower Chaussee
25, 12489 Berlin, Germany
\& ICMPA-UNESCO Chair, 072BP50, Cotonou, \& Ecole Normale Superieure, B.P 72, Natitingou, Benin}
\email{avohouremicocou@yahoo.fr}

\author{Joseph Ben Geloun}
\address[J.B.G.]{
Laboratoire d'Informatique de Paris Nord, UMR CNRS 7030
Universit\'{e} Sorbonne Paris Nord, 99, avenue J.-B. Clement, 93430 Villetaneuse, France, and
International Chair in Mathematical Physics and Applications,
ICMPA-UNESCO Chair, 072BP50, Cotonou, Rep. of Benin}
\email{bengeloun@lipn.univ-paris13.fr}

\author{Mahouton N. Hounkonnou}
\address[M.N.H.]{International Chair in Mathematical Physics and Applications,
ICMPA-UNESCO Chair, 072BP50, Cotonou, Rep. of Benin}
\email{mn.hounkonnou@cipma.net}

\begin{document}

\maketitle 

\begin{abstract}
In this paper, we analyze the Bollob\'as and Riordan polynomial $\cR$ for ribbon graphs with  half-ribbons introduced in [Combinatorics, Probability and Computing 31, 507-549, 2022]. We prove the universality property of a multivariate version of $\cR$ whereas $\cR$ itself turns out
to be universal for a subclass of ribbon graphs with 
half-ribbons. 
We also show that  $\cR$  can be defined on some equivalence classes of 
ribbon graphs involving half-ribbons moves and that the new polynomial is universal on these classes.\\

\end{abstract}

\tableofcontents

\section{Introduction}
The Bollob\'as-Riordan (BR) graph polynomial \cite{bollo} is a polynomial in four variables which extends the Tutte polynomial \cite{tutte,riv} from simple graphs to graphs with additional structures such as ribbon graphs (such graphs arise as neighbourhoods of graphs embedded into surfaces). Both polynomials satisfy a contraction/deletion recurrence rule defined on the associated graphs and, furthermore, are universal polynomial invariants. The universality property of these invariants means that any invariant of graphs satisfying the same relations of contraction and deletion can be calculated from those. Universality can be also of great use, for example, in statistical mechanics \cite{sok}
and quantum field theory \cite{Duchamp:2013pha,Duchamp:2013joa,Tanasa:2012pm}.

The BR polynomial is defined on signed ribbon graphs which are ribbon graphs whose edges are marked either by $+1$ or by $-1$. The signs of the edges play an important role in the orientability of the ribbon graphs. Signed ribbon graphs and their polynomial invariants are still under investigations \cite{vigne,joan,Tanasa:2010me,avo}. For example, in \cite{joan} the authors provide a ``recipe theorem'' for the BR polynomial  very close to the universality property. The proof of the universality of the BR polynomial is mainly based on the fact that the BR polynomial satisfies a contraction/deletion relation. However the proof of that claim relies on several other ingredients. Chord diagrams associated with bouquets and canonical diagrams found from these chord diagrams after a sequence of operations called rotations and twists about chords are extremely useful to establish that fact.

 Let us discuss in greater detail the polynomial on a new class of ribbon graphs introduced in \cite{rca} called ribbon graphs with  half-ribbons. A  half-ribbon (HR)  is simply a ribbon edge incident to a unique vertex without forming a loop. The presence of HRs in a ribbon graph has several interesting combinatorial properties as shown in \cite{rca}. HRs also allow to introduce a new and intuitive enough  operation which is the cut of an edge which differs from the usual edge deletion. The authors 
of the above work describe the implications that HRs have on the BR polynomial. One notes that in the polynomial worked out therein, the orientability of the ribbons is not taken into account. Since this new invariant satisfies a contraction/cut recurrence relation (replacing in this setting the usual contraction/deletion rule), one may wonder if this invariant is universal or not. Answering this question is the purpose of this paper.

We find in this paper an extension of the polynomial found in \cite{rca} by adding now a variable for the orientability of the ribbon graphs. We call it $\cR$. In the presence of this new variable, the contraction/cut rule still holds for $\cR$. We infer a
multivariate polynomial $\cQ$ that reduces to $\cR$ and also obeys a contraction/cut relation. 
 We then prove a main result (Theorem \ref{theo:univ2}) which is the universality property for $\cQ$ on ribbon graphs with HRs.  The method used to prove this is close to that given in \cite{bollo} but it is however specific due to the presence of HRs. As a corollary, the polynomial invariant $\cR$ is proved universal on 
 a subclass of ribbon graphs with HRs. 
 We then reveal the existence of another polynomial invariant defined over classes of ribbon graphs with HRs related to a new operation called HR moves. Theorem \ref{theo:univ} establishes the universality of that new polynomial which is a second main result of this paper.

The rest of this paper is organized as follows. In section \ref{sect:BRpoly}, we give an overview of the BR polynomial and its universality property. Section \ref{sect:bropen} recalls some results on the BR polynomial $\cR$ for ribbon graphs with HRs and introduces two multivariate versions 
$\widetilde{\cQ}$ and $\cQ$ of that invariant. The next section \ref{sect:univ} delivers our main result which is the proof of the universality theorem of $\cQ$. 
We finally define a polynomial invariant on classes of ribbon graphs related by moves of HRs and prove its universality property in section \ref{sect:moves}. 

\section{Overview of the Bollob\'as-Riordan polynomial and its universality property}
\label{sect:BRpoly}

In this section, we give an overview of the BR polynomial for ribbon graphs and mainly  focus on its  universality theorem introduced in \cite{bollo}. There are several ingredients in the proof of this theorem which will be useful for our subsequent developments and, thus, are worth to be reviewed as well.
\begin{definition}[Ribbon graphs \cite{bollo,krf}]\label{def:ribbongraph}
	A ribbon graph $\cG$ is  a (not necessarily orientable) surface with boundary represented as the union of two 
sets of closed topological discs called vertices $\cV$ and edges $\cE.$ These sets satisfy the following properties:

$\bullet$ Vertices and edges intersect in  disjoint line segments,

$\bullet$ each such line segment lies on the boundary of precisely one vertex and one edge,

$\bullet$ every edge contains exactly two such line segments.
\end{definition}

The isomorphism class of ribbon graphs is defined as follows \cite{bollo,bollo2}:
first identify a ribbon graph with its signed rotation system. 
Two rotation systems are called equivalent (henceforth their respective ribbon graphs) if there is a sequence of vertex flips and graph isomorphims that transform one into the other.

\begin{figure}[h]
\centering
\begin{minipage}[t]{.8\textwidth}
\centering
\includegraphics[angle=0, width=7cm, height=1cm]{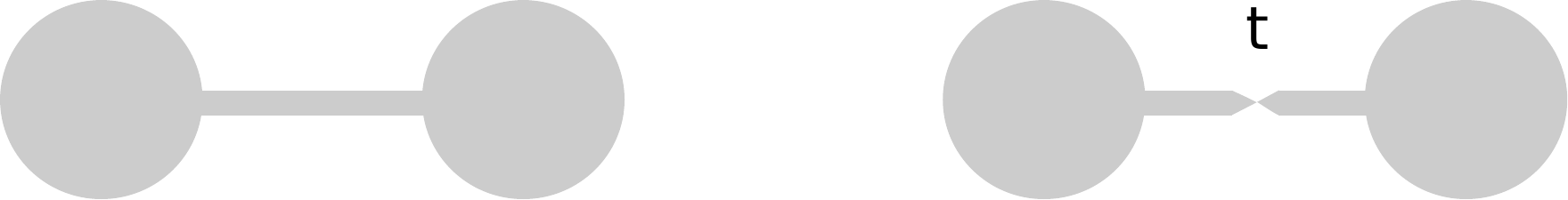}
\caption{ {\small Untwisted (left) and twisted (right) edge notations. }}
\label{fig:edges}
\end{minipage}
\end{figure}

There are three kinds of edges that can be identified in a ribbon graph.
An edge $e$ of a ribbon graph $\cG$ is called a bridge in $\cG$ if its removal disconnects a component of $\cG$. The edge $e$ is a self-loop in $\cG$ if the two ends of $e$ are incident to the same vertex $v$ of $\cG$ and $e$ is a regular edge of $\cG$ if it is neither a bridge nor a self-loop.  Ribbon edges can be  twisted as well (see Figure \ref{fig:edges}). We say that a self-loop $e$ at a vertex $v$ of a ribbon graph $\cG$ is twisted if $v\cup e$ forms a M\"obius band as opposed to an annulus (an untwisted self-loop). A self-loop $e$ is trivial if there is no cycle in $\cG$ which can be contracted to form a loop $f$ interlaced with $e$. Introducing twisted edges has some  consequences on the orientation of the ribbon graph.

In addition, there are other topological notions in a ribbon graph that we now describe. 

\begin{definition}[Faces and orientation \cite{bollo}]
A face is a component of the boundary of $\cG$ considered as a geometric ribbon graph and hence as a surface with boundary. 
\end{definition}
If $\cG$ is regarded as the neighborhood of a graph embedded into a surface, the set of faces is the set of faces of the embedding. 
A ribbon graph is  denoted by $\cG(\cV,\cE)$.

\begin{definition}[Deletion and contraction \cite{bollo}] 
\label{def:cdelrib}

Let $\cG$ be a ribbon graph and $e$ one of its edges. 

$\bullet$ We call $\cG-e$ the ribbon graph obtained from 
$\cG$ by deleting $e$ and keeping the end vertices 
as closed discs. 

$\bullet$ If $e$ is not a self-loop, the graph $\cG/e$  obtained by contracting $e$ is defined from $\cG$ by deleting $e$ and 
identifying its end vertices $v_{1,2}$ into a new vertex which possesses
all edges in the same cyclic order as they appeared in $v_{1,2}$.

$\bullet$ If $e$ is a trivial twisted self-loop, 
contraction is deletion: $\cG-e = \cG/e$. 
The contraction of a trivial untwisted self-loop  $e$
is the deletion of the self-loop and the addition of a vertex $v_0$ forming a new connected component to the graph  $\cG-e$. We write $\cG/e = (\cG-e)\sqcup \{v_0\}$.
\end{definition}

The contraction of a non self-loop may be also restated as follows: 
$\cG/e$ is defined from $\cG$ by  identifying a new vertex as $v_1\cup v_2\cup e$.
We recall that the contraction of a (twisted or untwisted) self-loop $e$
in $\cG$ coincides with an edge deletion in the graph dual of $\cG$.

A spanning subgraph $A$ of a ribbon graph $\cG(\cV,\cE)$ is a ribbon graph with set of vertices $\cV(A)=\cV$ and set of edges $\cE(A)\subseteq \cE$. We denote it as $A\sset \cG$.

\begin{definition}[Ribbon graph polynomial \cite{bollo}]
\label{def:bol}
Let $\cG$ be a ribbon graph. We define the ribbon graph polynomial of $\cG$ to be
\beq
R_{\cG}(X,Y,Z,W) =\sum_{A\sset \cG} (X-1)^{\rk(\cG)-\rk(A)}(Y-1)^{n(A)}Z^{k(A)-F(A)+n(A)}W^{t(A)}
\label{brpoly}
\eeq
considered as an element of the quotient of $\mathbb{Z}[X,Y,Z,W]$ by the ideal generated by $W^2-W$ 
and where $\rk(A),\,$ $n(A),\,$ $k(A),\,$ $F(A)\,$ and $t(A)$ are, respectively, the rank, the nullity, the number of connected components, 
the number of faces and the parameter which characterizes the orientability  of $A$ as a surface. If $A$ is orientable, then $t(A)=0$,  
otherwise, $t(A)=1$. By definition, $\rk(A)=|\cV|-k(A)$ and $n(A)=|\cE(A)|-\rk(A)$.
\end{definition}
In the following, we use the variable $(Y-1)$ for parameterizing the 
nullity of the subgraphs. This convention differs from the one in \cite{bollo} which rather uses $Y$. 
From a simple change of variable at any moment ($Y\rightarrow Y+1$), one can recover the convention used therein. 
Moreover, putting $W=1=Z$, one recovers the Tutte polynomial for $\cG$ seen as a simple graph. 
After introducing terminal forms, the choice $(Y-1)$ will be discussed.
We will often refer to the ribbon graph polynomial as the BR polynomial. Moreover, we use
interchangeably $\cR_{\cG}(X,Y,Z,W)$
and $\cR(\cG; X,Y,Z,W)$.

The BR polynomial obeys a contraction and deletion rule. 

\begin{theorem}[Contraction and deletion \cite{bollo}]
Let $\cG$ be a ribbon graph. If $e$ is a regular edge, then
\beq\label{regu}
R_{\cG}=R_{\cG/e}+ R_{\cG-e}\,;
\eeq 
for a bridge $e$ of $\cG$, one has
\beq
R_{\cG}=X \, R_{\cG/e}\,;
\label{cbrid}
\eeq 
for a trivial untwisted self-loop $e$,
\beq
R_{\cG}=Y \,R_{\cG-e}\,;
\label{csel}
\eeq
and for a trivial twisted self-loop $e$, the following holds
\beq
R_{\cG}=(1+(Y-1)ZW) \,R_{\cG-e}\,.
\label{ctsel}
\eeq
\end{theorem}

The relations \eqref{cbrid}--\eqref{ctsel} are useful for the evaluation of the terminal forms (ribbon graphs which only possess edges which are not regular). For a graph $\cG$ with only $n$ bridges, $m$ untwisted trivial self-loops and $p$ twisted trivial self-loops, the polynomial of $\cG$ is $X^nY^m(1+(Y-1)ZW)^p$. Note that, in \cite{avo}, the list of terminal forms has been further extended to specific one-vertex graphs called flowers so that one can complete the above with other contributions. 

Let us discuss in more details the universality of the BR polynomial for ribbon graphs \cite{bollo}. It is shown that the polynomial $R$ is the universal invariant for connected ribbon graphs satisfying  \eqref{regu} and \eqref{cbrid} and any other invariant satisfying the same relations can be calculated from $R$. First, one must understand that the knowledge of $R$ can be reduced to one-vertex ribbon graphs also simply called "bouquets".

Specifically, we obtain a bouquet  after a contraction of a spanning tree in a connected ribbon graph $\cG$. To achieve the proof of the universality of their polynomial, Bollob\'as and Riordan used another representation of  bouquets called ``signed chord diagrams'' (chord diagrams are also related to Vassiliev invariants \cite{bar,bir}). A chord diagram $D$ is a construction related to a bouquet $\cG$ such that if $\cG$ has $n$ edges, $D$ is constructed by putting on a circle $2n$ distinct points paired off by $n$ chords.  In the case of a ribbon graph with twisted and untwisted edges, $D$ is called a signed chord diagram,  if we put an assignment of sign ``$\tw$'' or ``$\unt$'' to each chord according to the fact that this chord corresponds to a twisted or negative edge or untwisted or positive edge, respectively.

 We shall write $n(D)$ for the number of chords of $D$ which is also the nullity of $\cG$ (each chord corresponds to an edge in a bouquet or a cycle generator). Using the ``doubling operation'' which consists in replacing each chord of $D$ by two edges joining the parts of the circle on each side of each end of the chord as shown in Figure \ref{fig:doubling}, $F(D)$ denotes the number of components of the resulting figure. We have $F(D)=F(\cG)$ and $t(D)$ stands for $t(\cG)$ which is equal to $0$ if all chords of $D$ have a positive sign (or untwisted) and $1$ otherwise.
%%%%%%%%%%%%%%%%%%%%%%%%%%%%%
\begin{figure}
    \label{fig:doubling}
\begin{tikzpicture}
  \draw (-1,0) -- (1,0);
  \draw (-1,-0.2) -- (1,-0.2);
  %\draw (0,-1.5) -- (0,1.5);
  \draw[dashed] (-1,0) .. controls (-1,0.555) and (-0.555,1) .. (0,1)
               .. controls (0.555,1) and (1,0.555) .. (1,0);
               
  \draw[dashed] (-1,-0.2) .. controls (-1,-0.755) and (-0.555,-1.2) .. (0,-1.2)
               .. controls (0.555,-1.2) and (1,-0.755) .. (1,-0.2); 
%%%%%%%%%%%%%%%%%%%%%               
  \draw (3,0) -- (5,-0.2);
  \draw (3,-0.2) -- (5,0);
  %\draw (0,-1.5) -- (0,1.5);
  \draw[dashed] (3,0) .. controls (3,0.555) and (3.445,1) .. (4,1)
               .. controls (4.555,1) and (5,0.555) .. (5,0);
               
  \draw[dashed] (3,-0.2) .. controls (3,-0.755) and (3.445,-1.2) .. (4,-1.2)
               .. controls (4.555,-1.2) and (5,-0.755) .. (5,-0.2); 
\end{tikzpicture}
    \caption{A doubled positive chord  (left) and a doubled negative chord  (right).}
\end{figure}

%%%%%%%%%%%%%%%%%%%%%%%
A subdiagram of a signed chord diagram $D$ is a signed chord diagram $D'$ obtained from $D$ by deleting a subset of chords of $D$. For a bouquet $\cG$, looked at as a signed chord diagram $D$, the BR polynomial summation is defined over the spanning subdiagrams $D'\sset D$ as:
\bea
R(D;X,Y,Z,W)=\sum_{D'\sset D}(Y-1)^{n(D')}Z^{1-F(D')+n(D')}W^{t(D')}.
\label{sumdi}
\eea
Later this summation is written as:
\bea
R(D;X,Y,Z,W)=\sum_{i,j,k}R_{ijk}(D; X)(Y-1)^iZ^{j}W^k,
\label{sumdic}
\eea
where $R_{ijk}(D;X)$, the coefficient of $(Y-1)^iZ^{j}W^k$ in \eqref{sumdic},  counts the number of subdiagrams $D'\sset D$ which have $i$ chords, $j=1-F(D')+n(D')$ and $k=t(D')$ in \eqref{sumdi}.
It is obvious that the
above expression finds an extension to any ribbon graph $\cG$.
In such a case, 
$R_{ijk}(D;X)$ 
becomes a sum of monomials $(X-1)^{\rk(\cG)-\rk(A)}$
for particular subgraphs $A\sset \cG$ with properties constrained by $i,j,k$.

Let $\mathfrak{G}$ be the set of isomorphism classes of connected ribbon graphs \cite{bollo}. The theorem of universality is given by the following statement:

\begin{theorem}[Universality of Bollob\'as-Riordan polynomial \cite{bollo}]
\label{theo:univ1}
Let $\mathfrak{R}$ be a commutative ring, $x$ an element of $\mathfrak{R}$, and $\phi$ a map from $\mathfrak{G}$ to $\mathfrak{R}$ satisfying
\bea \label{system1}
\phi(\cG) = \left\{\begin{array}{ll} 
\phi(\cG - e) + \phi(\cG/e) & {\text{if e is  regular}},\\\\
x \, \phi(\cG/e) & {\text{if e is a bridge}}.
\end{array} \right. 
\eea
\medskip 
Then there are elements $\lambda_{ijk}\in \mathfrak{R}$, $i\geq 0$, $0\leq j\leq i$, $0\leq k \leq 1$, such that 
\beq \label{univsum1}
\phi(\cG) = \sum_{i,j,k} \lambda_{ijk} R_{ijk}(\cG;x). 
\eeq
\end{theorem}

The main point of the universality theorem is the determination of the $\lambda_{ijk}$. 
 The coefficients $\lambda_{ijk}$ are   determined by the evaluation of $\phi$ on the so-called ``canonical diagrams''. 
\begin{figure}[h]
 \centering
     \begin{minipage}[t]{.8\textwidth}
      \centering
\includegraphics[angle=0, width=1.5cm, height=1.5cm]{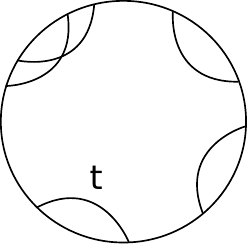}
\caption{ {\small The canonical chord diagram $D_{5,1,1}$. }}
\label{fig:canodi}
\end{minipage}
\end{figure}
For a bouquet $\cG$ seen as a chord diagram $D$, a sequence of rotations and twists about chords \cite{bollo} in $D$ provides a simple diagram called canonical. 
Given canonical diagrams $D_{i,j,k}$, consisting of $i-2j-k$ positive chords intersecting no other
chords, $j$ pairs of intersecting positive chords, and $k$ negative chords  $0\leq k\leq2$, intersecting no other chords (see an example in Figure \ref{fig:canodi}), then  $\lambda_{ijk}$ is equal to some $\phi(D_{i,j',k'})$. This is proved by a recurrence relation on the number of chords $i$, given the initial value $\lambda_{000}$ for the value of $\phi$ on a bare vertex. The same result holds  for any connected ribbon graph using the relations \eqref{system1}. The case of several connected components can be simply inferred from this point because the polynomial is multiplicative over disjoint union.
 
\section{The Bollob\'as-Riordan polynomial for
 ribbon graphs with HRs}
\label{sect:bropen}
This section introduces a polynomial invariant for ribbon graphs with HRs which is a notion studied in \cite{krf}.  
The polynomial which will be discussed extends the invariant found in \cite{rca} by adding a  variable that takes into account the orientability of the ribbon graph. We introduce
some multivariate variants of that polynomial. 
The question of the universality
of such polynomials 
is then asked.  

We first recall some definitions.

\begin{figure}[h]
 \centering
     \begin{minipage}[t]{.8\textwidth}
      \centering
\includegraphics[angle=0, width=3cm, height=1.5cm]{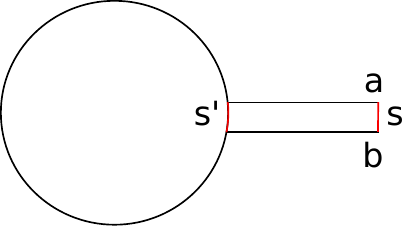}
\caption{ {\small A HR with two end segments (in red): $s'$
touching the vertex and $s$ external; the ends $a$ and $b$ of $s$
are the external points. }}
\label{fig:HR}
\end{minipage}
\end{figure}

\begin{definition}[Half-ribbon and external points \cite{rca}]
\label{ribHR}
A half-ribbon or half-edge  is a rectangle  incident to a unique vertex of a ribbon graph by
a unique line segment $s$ on the boundary, i.e. without forming a loop. The segment parallel to $s$ called the external segment. The end points of any external segment are called external points
of the HR. The two boundary segments of a ribbon edge or of a HR that are neither external nor
incident to a vertex are called strands. A HR is always oriented consistently with the vertex it
intersects. 
   (See Figure \ref{fig:HR}.)
\end{definition}

\begin{figure}[h]
 \centering
\begin{minipage}[t]{.8\textwidth}
\centering
\includegraphics[angle=0, width=6cm, height=1cm]{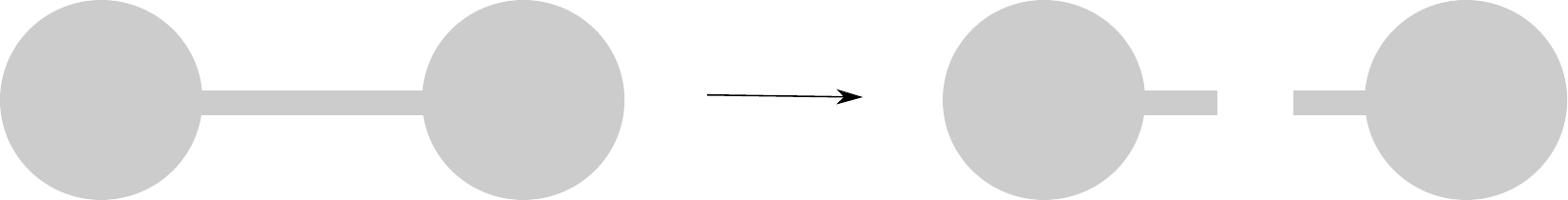}
\caption{{\small Cutting a ribbon edge.}}
\label{fig:cut}
\end{minipage}
\end{figure}

\begin{definition}[Cut of a ribbon edge \cite{rca}]
\label{cutriedge}
Let $\cG$ be a ribbon graph and let $e$ be a ribbon edge of $\cG$.
The \emph{cut graph} $\cG \vee e$ is obtained from $\cG$ by 
deleting $e$ and attaching two HRs
at the same line segments where $e$ was incident
 to the end vertices, one at each of the end vertices of $e$. If $e$ is a loop, the two HRs are on the same  vertex.    (See an illustration in Figure \ref{fig:cut}.)
\end{definition}

\begin{figure}[h]
 \centering
     \begin{minipage}[t]{.8\textwidth}
      \centering
\includegraphics[angle=0, width=6.5cm, height=2cm]{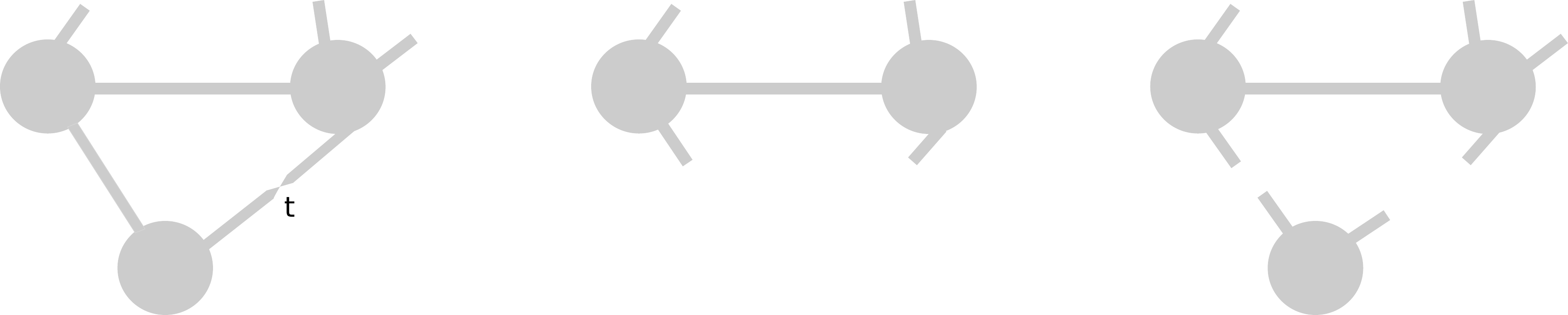}
\vspace{0.3cm}
\caption{{\small A ribbon graph with HRs $\cG_{\mf^0}$ and
a spanning c-subgraph $A$.}}
\label{fig:ribgraph}
\end{minipage}
\put(-225,-10){$\cG$}
%\put(-175,-10){$A$}
\put(-107,-10){$A$}
\end{figure}

The definition of a ribbon graph with HRs may be introduced at this stage.
\begin{definition}[Ribbon graph with HRs and spanning c-subgraph \cite{rca}]
\label{def:herg}
 A ribbon graph with HRs $\cG(\cV,\cE,\mf^0)$ (or simply $\cG_{\mf^0}$) 
is a ribbon graph $\cG(\cV,\cE)$ (or shortly $\cG$) with a set $\mf^0$ of HRs
 such that each HR is attached to a unique vertex as in Definition \ref{ribHR}, and the segments where the HRs are attached 
are disjoint from each other and from the segments where any ribbon edges are attached. 
The ribbon graph $\cG$  is called the underlying ribbon graph of $\cG_{\mf^0}$.

$\bullet$ A spanning c-subgraph
$A$ of $\cG_{\mf^0}$ is formed by cutting some subset of the ribbon edges of $\cG_{\mf^0}$. We denote again the spanning c-subgraph inclusion as $A \sset \cG_{\mf^0}$. 
(See $A$ in Figure \ref{fig:ribgraph}.)
\end{definition} 

Note that a ribbon graph is a ribbon graph with HRs with $\mf^0=\emptyset$.
The isomorphism class of ribbon graph with HRs is much inspired
from the isomorphism class of ribbon graphs. Consider two ribbon graphs with HRs $\cG_{\mf^0}$ and $\mathcal{H}_{\mf'^0}$. 
We say that $\cG_{\mf^0}$ is isomorphic to 
$\mathcal{H}_{\mf'^0}$, 
if their underlying ribbon graphs $\cG$ and $\mathcal{H}$ are isomorphic, and their sets of HRs are of same cardinality and obeys the same incidence relation 
with the same cyclic ordering onto vertices. 

The  notion that we will extensively use is the one of spanning c-subgraph. We can simply explain that notion in the following way: 
Take a subset of edges of a given graph, cut them all. Consider the spanning subgraph then formed by the resulting graph. 
The set of HRs of this subgraph contains both the set of HRs of the initial graph ($\mf^{0}$) plus an additional set induced by the cut of the edges. 

Note that cutting an edge of a graph modifies the boundary faces of this graph. There are new boundary faces following the contour of the HRs.  
However, combinatorially, we distinguish this new type of faces and the initial ones which follow the boundary of well-formed edges.

\begin{figure}[h]
 \centering
     \begin{minipage}[t]{.8\textwidth}
      \centering
\includegraphics[angle=0, width=4cm, height=2cm]{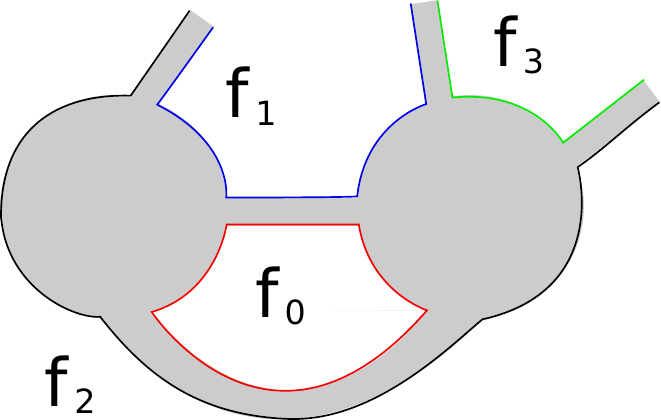}
\vspace{0.3cm}
\caption{ {\small A ribbon graph with set 
of internal faces 
$\cF_{\inter}=\{f_0\}$, and set of external 
faces $\cF_{\ext}=\{f_1,f_2,f_3\}$. }}
\label{fig:faces}
\end{minipage}
\end{figure}

\begin{definition}[Closed, open faces]
\label{faces}
Let $\cG_{\mf^0}$ be a ribbon graph with HRs. 

$\bullet$ A \emph{closed face} is a boundary  component of $\cG_{\mf^0}$  which never passes through any external segment of a HR.  The set of closed faces is denoted $\cF_{\inter}$. (See the closed face $f_{0}$ in 
Figure \ref{fig:faces}.)

$\bullet$ An \emph{open  face}  is a boundary arc
leaving an external point of some HR rejoining another external point
 without passing through any external segment of a HR.  
 The set of open faces is denoted  $\cF_{\ext}$. 
  (Examples of open faces are provided in Figure \ref{fig:faces}.)

$\bullet$ The set of faces $\cF$ of a ribbon graph with HRs  is defined by 
$\cF_{\inter} \cup \cF_{\ext}$.
\end{definition}

Open and closed faces are illustrated in Figure \ref{fig:faces}.

\begin{figure}[h]
 \centering
     \begin{minipage}[t]{.8\textwidth}
      \centering
\includegraphics[angle=0, width=2cm, height=1cm]{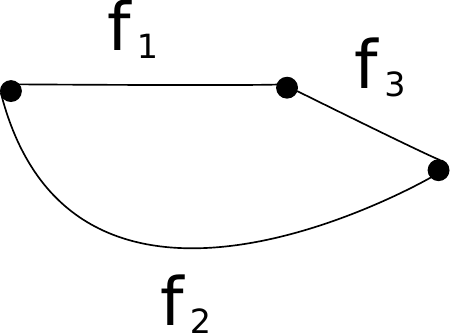}
\vspace{0.3cm}
\caption{ {\small The boundary graph associated with 
the ribbon graph in Figure \ref{fig:faces}. }}
\label{fig:boundary}
\end{minipage}
\end{figure}

\begin{definition}[Boundary graph \cite{Gurau:2009tz}]
\label{bnd}
The \emph{boundary graph} $\bG$  of a ribbon graph with HRs $\cG_{\mf^0}$ is an abstract graph $\bG(\bV,\bE)$ 
such that $\bV$ is in one-to-one correspondence with $\mf^0$,
and $\bE$ is in one-to-one correspondence with $\cF_{\ext}$. 
Consider an edge $e$ of $\bE$, its corresponding open face $f_e \in\cF_{\ext} $, a vertex $v$, and its corresponding 
HR $h_v$. The edge $e$ is incident to $v$ if and only if  $f_e$ has one end-point in $h_v$,  and,  if both end-points of $f_e$ are in $h_v$, then $e$ is a loop.   (The boundary of the graph given in Figure \ref{fig:faces}
is provided in  Figure \ref{fig:boundary}.)
\end{definition} 

The notions of edge contraction and deletion for ribbon graphs with HRs keep their meaning as in Definition \ref{def:cdelrib}. 
We are in position to identify a new polynomial invariant. To alleviate our notation, 
from this point onwards, 
we will denote a ribbon graph with HRs $\cG_{\mf^0}$ simply 
as $\cG$ as there will be no confusion given the fact
that we will always work with ribbon graph with HRs.
We keep the notation 
$\mathfrak{G}$ for the set of isomorphism classes of connected ribbon graphs with HRs.

\begin{definition}[Polynomial for ribbon graphs with HRs]
Let $\cG(\cV,\cE,\mf^0)$ be a ribbon graph with HRs.  
We define the polynomial of $\cG$ to be
\beq
\cR_{\cG}(X,Y,Z,S,W,T)
=\sum_{A\sset \cG} (X-1)^{\rk(\cG)-\rk(A)}(Y-1)^{n(A)}
Z^{k(A)-F_{\inter}(A)+n(A)} \, S^{C_\partial(A)} \, W^{t(A)}T^{f(A)},
\label{brfla}
\eeq
considered as an element of the quotient of $\mathbb{Z}[X,Y,Z,S,W,T]$ by the ideal generated by $W^2-W$, 
where $C_\partial(A)= |\bC(A)|$ is the number of connected components of the boundary graph of $A$, 
$F_{\inter}(A)= |\cF_{\inter}(A)|$ and $f(A)$ the number of HRs of $A$.
\end{definition}

The polynomial $\cR$ \eqref{brfla} generalizes the BR polynomial 
$R$ \eqref{brpoly}.  
One way to recover $R$ is by replacing the sum over c-subgraphs 
by a sum over subgraphs (that casts the powers of $S$ to 0), then putting $T=1$. 
Another remark follows: as the number of HRs on a spanning c-subgraph $A$
 can be written as $f(A)= |\mf^0|+2(|\cE| - |\cE(A)|)$, 
 we can always factor 
 $T^{f(\cG)}$ from the polynomial. Hence, 
 $\cR_{\cG}/T^{|\mf^0|}$
 might be also
 considered as an interesting reduced polynomial.
 
 After performing the change of variable $S\to Z^{-1}$,
we are led to another extension of the BR polynomial for ribbon graphs with HRs. 
We will refer the second polynomial to as $\cR'$. 
In symbols, for a ribbon graph with HRs $\cG$, we write
\beq\label{rpr}
\cR_{\cG}(X,Y,Z,Z^{-1},W,T) = \cR'_{\cG}(X,Y,Z,W,T) \,,
\eeq
where $\cR$ is given by Definition \ref{brfla} .

Graph operations such as the disjoint union and the one-point-join ($\cG_1 \sqcup \cG_2$ and $\cG_1 \cdot_{v_1,v_2} \cG_2$, respectively) \cite{tutte} 
extend to ribbon graphs \cite{bollo} and to ribbon graphs with HRs \cite{rca}. The product $\cG_1 \cdot_{v_1,v_2} \cG_2$ 
at the vertex resulting from merging $v_1$ and $v_2$ on an arc of each
of these which does not contain any ribbon edges or HRs
(in the sense of the second point of Definition \ref{def:cdelrib}) 
respects the cyclic order of all edges and HRs on the previous vertices $v_1$ and $v_2$. The fact that 
$R_{\cG_1 \sqcup \cG_2}= R_{\cG_1} R_{\cG_2} = R_{\cG_1 \cdot_{v_1,v_2} \cG_2}$ holds
for ribbon graphs without HRs \cite{bollo} can be 
extended  to ribbon graphs with HRs under particular conditions. The following proposition holds.
\begin{proposition}[Operations on BR polynomials \cite{rca}] 
\label{opBRfla}
Let $\cG_1$ and $\cG_2$ be two disjoint ribbon graphs with HRs, then 
\bea 
\cR_{\cG_1 \sqcup \cG_2}&=&\cR_{\cG_1}
 \cR_{\cG_2},\quad \cR'_{\cG_1 \sqcup \cG_2}= \cR'_{\cG_1}
 \cR'_{\cG_2}\,,
\label{brcups1}\\
\cR'_{\cG_1 \cdot_{v_1,v_2} \cG_2} 
&=&\cR'_{\cG_1} \cR'_{\cG_2} \,,
\label{brcups2}
\eea
for any disjoint vertices $v_{1,2}$ in $\cG_{1,2}$, respectively. 
\end{proposition}
\proof  The proof of Proposition  \ref{opBRfla} corresponds to that of Proposition 5 in \cite{rca} 
where the sole additional fact concerns the variable $W$ associated with the orientability. 
This can be simply achieved by adding the fact that $W^2=W$ in the proof of Proposition 5 in \cite{rca}. 
\qed

\begin{theorem}[Contraction and cut on BR polynomial]
\label{theo:BRext}
Let $\cG(\cV,\cE,\mf^0)$ be a ribbon graph with HRs. Then,
for a regular edge $e$, 
\beq
\cR_{\cG}=\cR_{ \cG\vee e} +\cR_{\cG/e}\,;
\label{retcondel}
\eeq 
for a bridge $e$, we have 
\beq
\cR_{\cG} =(X-1)\cR_{\cG \vee e}+ \cR_{\cG/e}\,;
\label{retbri}
\eeq 
 for a trivial twisted self-loop $e$, the following holds 
\beq
\cR_{\cG}=\cR_{\cG \vee e} + (Y-1)ZW \,\cR_{\cG/e}\,;
\label{retsel2}
\eeq 
whereas for a trivial untwisted self-loop $e$, we have
\beq
\cR_{\cG}=\cR_{\cG \vee e} + (Y-1) \cR_{\cG/e}\,.
\label{retsel1}
\eeq 
\end{theorem}

\proof  This can be proved in the same lines of Theorem 3 in \cite{rca} where the new point \eqref{retsel2} 
associated with the orientability can be recovered from \cite{bollo}. 
\qed

\begin{corollary}[Contraction and cut on BR polynomial $\cR'$]
\label{coro:rprim}
Let $\cG(\cV,\cE,\mf^0)$ be a ribbon graph with HRs. Then,
for a regular edge $e$, 
\beq
\cR'_{\cG}=\cR'_{ \cG\vee e} +\cR'_{\cG/e}\,,
\qquad 
\cR'_{\cG \vee e} = T^2 \, \cR'_{\cG -e}\,;
\label{retcondelp}
\eeq 
for a bridge $e$, we have 
$\cR'_{\cG/e} = \cR'_{\cG-e} = T^{-2} \, \cR'_{\cG \vee e}$
\beq
\cR'_{\cG} =[(X-1)T^2 +1]  \, \cR'_{\cG/e}\,;
\label{retbrip}
\eeq 
for a trivial twisted self-loop, $\cR'_{\cG-e}=T^{-2} \, \cR'_{\cG \vee e}$ and 
\beq
\cR'_{\cG}=[T^2 + (Y-1)ZW] \,\cR'_{\cG-e}\,,
\label{retselp2}
\eeq 
whereas for a trivial untwisted self-loop, we have
$\cR'_{\cG-e}=T^{-2} \, \cR'_{\cG \vee e}$
and 
\beq
\cR'_{\cG}=[T^2+ (Y-1)] \, \cR'_{\cG-e}\,.
\label{retselp1}
\eeq 
\end{corollary}

\proof 
The corollary is immediate from Theorem \ref{theo:BRext} and Corollary 1 in \cite{rca}. 
The new relation \eqref{retselp2} can be achieved using a similar identity in Theorem 1 in \cite{bollo}.

\qed

The following statement holds.
\begin{proposition}
The polynomial $\cR'$ is universal
in the sense of Theorem \ref{theo:univ1}. 
\end{proposition}
\proof 
From  Corollary \ref{coro:rprim}, $\cR'$ satisfies the following relations:
\beq
\cR'_{\cG}=
\left\{
\begin{array}{ll} 
 T^2 \,\cR' (\cG - e)  + \cR' (\cG/e)  \qquad\quad  
{\mbox{if $e$ is neither a bridge nor a self-loop}}, \cr
[(X-1)T^2 +1] \, \cR'(\cG - e) \quad   {\mbox{\quad if $e$ is a bridge.}} 
\end{array}
\right. 
\label{univ}
\eeq
After a change of variables as:
\begin{eqnarray}
\left\{\begin{array}{ll}\tilde X=(X-1)T^2+1\\
 \tilde Y=Y-1+T^2\end{array}
\right.
\end{eqnarray}
and given the fact that, for a given ribbon graph $\cG(\cV,\cE,\mf^0)$
and $A \sset \cG$,  
\beq
f(A)=|\mf^0|+2(|\cE|-|\cE(A)|)\,,
\eeq
we get 
\beq
 \cR'_{\cG}(X,Y,Z,W,T)=T^{|\mf^0|}T^{2n(\cG)} \,R_{\cG}(\tilde X, \frac{\tilde Y}{T^2},Z,W), \,
\eeq
with $R$ the BR polynomial defined in \eqref{brpoly}. The above equation shows that the
reduced polynomial $\cR'$ is universal on the set of ribbon graphs, i.e. it defines 
a family of 
base polynomials $\{\cR'_{ijk}\}$ that 
plays the same role as the
family $\{R_{ijk}\}$ in Theorem \ref{theo:univ1}. 
\qed 
%
%{\cred 
%Another remark can be made: as $\cR'$ fulfills Corollary \ref{coro:rprim}, 
 %we may see it as a more natural invariant on the set of ribbon graphs with HRs on which it is also a universal invariant.
%}

\

\noindent{\bf Multivariate polynomials.}
There exist multivariate 
versions of $\cR$. We will concentrate on a general  multivariate polynomial and 
one of its reduction that turns out to be universal for the invariants satisfying
the contraction-cut rule on 
ribbon graphs with HRs. 

We introduce some notation. 
Let $\cG$ be a ribbon graph with HRs. Consider any 
$A\sset \cG$
and call $\mf(A)$ the set of HRs of $A$. We recall 
the notation $|\mf(A)|= f(A)$. 
%We  decompose 
%$\mf(A)= \mf^0\cup \mf^1(A)$ where
%$\mf^0\cap \mf^1(A)=\emptyset$; 
%  $\mf^1(A)$
% comes from the cut of a subset of edges in $\cE$, those that do not define $\cE(A)$.  
Let $\cG^0$  be the spanning c-subgraph of $\cG$ obtained by cutting all edges in $\cE$. 
Any HR of any spanning c-subgraph $A$ of $\cG$ must appear (once and only once) in $\mf(\cG^0)$. This also means 
$\mf(A) \subseteq \mf(\cG^0)$. 
The following polynomial 
requires at most  $|\mf(\cG^0)|$ variables 
for each of its monomials.

 For any c-subgraph $A\sset \cG$, $\cF_{\inter}(A) \subseteq 
 \cF_{\inter}(\cG)$, as a closed 
 face of $\cG$ could be either cut during the process of creating $A$ or kept in $A$. 
 Thus, we introduce a set of variables  $z_{\alpha}$ for $\alpha\in  \cF_{\inter}=\cF_{\inter}(\cG)$.

%{\color{green} For any c-subgraph $A\sset \cG$, consider the set $\bC(A)$ of connected components of the boundary of $A$. 
 %We recall $|\bC(A)|= C_\partial(A)$. 
%An element $\mathbf{c}$ of $\bC(A)$ can be labelled by $c \in \{1,\dots, C_\partial(A)\}$. 
%Therefore, we can denote it as $\mathbf{c}_c$. 
%However, we make a particular label choice in such a way that 
%the number $l_c$ of HRs that are incident to the element labeled by $c$ are ordered as $l_1 \ge l_2 \ge \dots \ge 
%l_{ C_\partial(A)}$. 
%Note that $\sum_{c=1}^{C_\partial(A)} l_c = f(A)$
%defines a partition denoted
%$[l_c]$ of the number of 
%HRs of $A$. For two elements of $\bC(A)$ 
%labeled by $c$ and $c'$, if  $l_c = l_{c'}$, then we consecutively order $c$ and $c'$ arbitrarily.  
%Consider $C_{\partial}^{\max}= \max_{A \sset \cG} C_\partial(A)$, 
% we introduce a set of variables $\{s_{c}\}_{c\in 
 %\{1,\dots, C_{\partial}^{\max} \}}$.
 
  %} 

For any c-subgraph $A\sset \cG$, consider the set $\bC(A)$ of connected components of the boundary of $A$. 
 We recall $|\bC(A)|= C_\partial(A)$. Let $\{ \zeta_i \}_{i = 0, \cdots, |f(A)|}$
 be a set of variables,  
 where each $\zeta_i$ records the presence of one boundary component with $i$ half-ribbons on it. 
 More explicitly, for a c-subgraph $A$ 
  and an element $\mathbf{c}_c$ of $\bC(A)$, we define $\mathbf{c}_c\cap \mf(A)$ to be the subset of HRs  of $A$ that are incident to $\mathbf{c}_c$.
To lighten the notation, we shortly write
$c\cap \mf(A)$ instead of $\mathbf{c}_c\cap \mf(A)$. 
We associate a variable $\zeta_{i=|c\cap \mf(A)|}$  with each $c \in \{1,\dots, C_\partial(A)\}$. Clearly,   $f(A)$ becomes partitioned by the $\zeta_i$ in $|C_\partial(A)|$ parts,
each of which associated with a connected component of the boundary graph: 
$f(A)=\sum_{c} |c\cap \mf(A)|$.

\begin{definition}[Multivariate polynomial for ribbon graphs with HRs]
Let $\cG(\cV,\cE,\mf^0)$ be a ribbon graph with HRs,
 $\{\beta_e\} = \{\beta_e\}_{e\in \cE}$ be a set of variables associated with the edges of $\cG$, 
$\{z_{\alpha}\} = \{z_{\alpha}\}_{\alpha \in \cF_{\inter}}$  be a set of variables associated with internal faces $\cG$.
Let $ \{\zeta_{i} \}_{ i\in \{1, \dots ,|\mf(\cG^0)|\}}$ be a set of variables recording $i$ HRs in any connected component of the boundary of any c-subgraphs of $\cG$. %Varying $A$, we denote by $\{\zeta_{c,f}\}$ the set all possible variables of that kind. 

We define the multivariate polynomial of $\cG $ 
 to be 
% \bea
%&&
%\widetilde\cQ_{\cG}(x,\{\beta_e\},\{z_\alpha\},s,w,\{\zeta_{i}\})
%=\crcr
%&&
%\sum_{A\sset \cG} x^{ k(A)}
%\left(\prod_{e\in \cE(A)} \beta_e\right)
%\left(
%\prod_{
%\alpha\in  \cF_{\inter}(A)}
%z_\alpha 
%\right)\,
%w^{t(A)}
%s^{C_\partial(A)}
%\left(
%\prod_{c=1}^{C_\partial(A)} s_{c} 
%\prod_{f\in \, c \cap \mf(A)}\zeta_{c,f}
%\right) .
%\label{brflam}
%\eea

 \bea
&&
\widetilde\cQ_{\cG}(x,\{\beta_e\},\{z_\alpha\},w, \{\zeta_i\})
=\crcr
&&
\sum_{A\sset \cG} x^{ k(A)}
\left(\prod_{e\in \cE(A)} \beta_e\right)
\left(
\prod_{
\alpha\in  \cF_{\inter}(A)}
z_\alpha 
\right)\,
w^{t(A)}\,
\left(\prod_{c=1}^{C_\partial(A)} \zeta_{|c \cap \mathfrak{f}(A)|}\right).
\label{brflam}
\eea

\end{definition}

The following statement holds 
\begin{theorem}
Let $\cG(\cV,\cE,\mf^0)$ be a ribbon graph with HRs. Then,
for a regular edge $e$, 
the multivariate polynomial obeys the recursion relation
\bea
&&
 \widetilde
\cQ_{\cG}(x,\{\beta_e\},\{z_\alpha\},w, \{\zeta_{i}\})
\\
&& = \widetilde
\cQ_{\cG \vee e }(x,\{\beta_{e'\ne e }\},\{z_\alpha\},w, \{\zeta_{i}\})
+
x \beta_e 
\widetilde
\cQ_{\cG/e}(x,\{\beta_{e'\ne e }\},\{z_\alpha\},w, \{\zeta_{i}\})
\nonumber
\eea
\end{theorem}
\proof 
The proof is straightforward as it follows the standard
decomposition of the set of c-subgraphs in 
those that contain $e$ and those
that do not.  
\qed

We have the following
reduction, setting all multiple
variables to some constants
\bea
&&
\widetilde\cQ_{\cG}(x,\{\beta_e= \beta\},\{z_\alpha=z\},w,\{\zeta_{i}=s\zeta^i\})=\cr\cr
&&
\sum_{A\sset \cG} x^{ k(A)}
\beta^{E(A)}
z^{F_{\inter}(A)} 
 w^{t(A)}\, s^{C_\partial(A)} \,
\zeta^{f(A)}\crcr
%&& = 
%(xz\beta^{-1})^{- %k(\cG)+k(\cG)}
%\sum_{A\sset \cG} x^{ k(A)}
%\beta^{-k(A) + |\cV|}
%\beta^{E(A)+k(A) - |\cV|}
%z^{F_{\inter}(A)} \, %s^{C_\partial(A)} \, w^{t(A)}
%\zeta^{f(A)} \crcr
&& = 
(zx\beta^{-1})^{k(\cG)}
\beta^{|\cV|}
\sum_{A\sset \cG} (xz\beta^{-1})^{ k(A)-k(\cG)}
(z\beta)^{n(A)}
z^{-n(A)-k(A)}  
z^{F_{\inter}(A)} \, s^{C_\partial(A)} \, w^{t(A)}
\zeta^{f(A)} \crcr
&&
= 
(zx\beta^{-1})^{k(\cG)}
\beta^{|\cV|}\, 
\cR_{\cG}(xz \beta^{-1}+1,z\beta+1,z^{-1},s,w,\zeta)
\eea
where $E(A) = |\cE(A)|$. Thus $\cR_{\cG}$ can be recovered from the multivariate polynomial $\widetilde\cQ_{\cG}$ after some change of variables. 

We will be interested in 
the intermediate reduced form 
\bea
&&
\frac{1}{((x-1)(y-1)z)^{k(\cG)}}
\frac{1}{((y-1)z)^{|\cV|-k(\cG)}} \times \crcr
&&
\widetilde \cQ_{\cG}((x-1)(y-1)z^2 ,
\{\beta_e=  (y -1)z \},
\{z_\alpha=z^{-1}\},w,\{\zeta_{i}\})
\cr\cr
&& 
= 
\frac{1}{((x-1)(y-1)z)^{k(\cG)}}
\frac{1}{((y-1)z)^{|\cV|-k(\cG)}} \times \crcr
&&
\sum_{A\sset \cG}  ((x-1)(y-1)z^2) ^{ k(A)}
( (y -1)z)^{E(A)}
z^{-F_{\inter}(A)}\,
w^{t(A)}
 \prod_{ c =1}^{C_\partial(A)} \zeta_{|c \cap \mathfrak{f}(A)|}
 \crcr
 && 
%= 
%\sum_{A\sset \cG} 
% (x-1) ^{ k(A)- k(\cG)}
%(y-1) ^{ k(A) + E(A) - |\cV|}
%(z^2) ^{ k(A)}
%( z)^{-|\cV| + E(A)}
%z^{-F_{\inter}(A)}\,
%w^{t(A)}
%s^{C_\partial(A)} 
%\prod_{ c =1)}^{C_\partial(A)} %\zeta_{c}^{l_c}
 %\crcr
%&&
= 
\sum_{A\sset \cG}  (x-1) ^{ k(A)- k(\cG)} (y -1)^{n(A)} 
z^{k(A)-F_{\inter}(A) + n(A)}\,
w^{t(A)} 
\prod_{ c =1}^{C_\partial(A)} \zeta_{|c \cap \mathfrak{f}(A)|}
\crcr
&&
 = \cQ_{\cG}(x,y,z,s,w,\{\zeta_{i}\})
\label{brflamulti}
\eea

Thus 
$ \cQ_{\cG}(x,y,z,s,w, \{\zeta_{i}\})$
defines a multivariate invariant 
with monomials that keep track of
the partition of the number of HRs 
on each connected component of the boundary graph of each c-subgraph. 
$\cR_{\cG}$ can be recovered 
from $\cQ_{\cG}$ setting
all $\zeta_i = s\zeta^i$.

\section{Main results: Universality theorems}
\label{sect:univ}

\subsection{Chord diagrams with HRs}

The main objective of this sub-section is the determination of a special class of diagrams called canonical  which turn out to be necessary for the proof of the universality of the polynomial in \eqref{brfla}. To succeed in this, we need to understand how the operations of rotation and twist about chords \cite{bollo} make sense on ``open'' chord diagrams or chord diagrams associated to one-vertex ribbon graphs with HRs, called bouquets with HRs. After defining open chord diagrams, we will focus on a two-vertex ribbon graph with HRs where the distinct ways of contracting the edges lead to some equivalent diagrams.

\begin{definition}[Chord diagrams]
$\bullet$ A HR on a chord diagram is a segment attached to a unique point on its circle. 

$\bullet$ An (open) chord diagram is a chord diagram in the sense of \cite{bollo} with the further data of the set of HRs. In the case where this set is empty, it becomes a chord diagram.

$\bullet$ A signed (open) chord diagram is an (open) chord diagram with an assignment of a sign ``\tw'' or ``\unt'' to each chord.
\end{definition}

We remark that in the previous definition of chord diagram $D$, if $D$ has $n$ chords and $l$ HRs, there are $2n+l$ distinct and marked points on the circle. 

If $\cG$ is a bouquet with HRs and $D$ the corresponding (open) signed chord diagram, the number  of chords $n(D)$ of $D$ is equal to the nullity of $\cG$ and we have $n(D)=e(\cG)=n(\cG)$. The doubling operation on $D$ consists of replacing each chord of $D$ by two edges joining the parts of the circle on each side of each end of the chord and each HR of $D$ by two parallel  segments,  each one
on each side of the HR. For each HR, 
consider the end points
of the two parallel segments
that are not on the circle. 
Insert a vertex of degree 2 
between these end points
and perform this insertion for 
each pair of parallel 
segments for each HR. 
We call the resulting
diagram $D'$, the pinched diagram of $D$. With this operation, the number of boundary components of  $D'$ is equal to $F_{\inter}(D')+C_\partial(D')$ where $F_{\inter}(D')=F_{\inter}(D) $ is the number of components which are closed and $C_\partial(D')$ is the number of remaining
boundary components. 
We then define $C_\partial(D) =C_\partial(D')$. 
We easily realize  that $C_\partial(D)$ 
is equal to the number 
of connected components of
the boundary graph associated
with $D$.

The ordinary operations on ribbon graphs simply translate to chord diagrams. In particular, the deletion or the cutting of chords and disjoint union or one-point-join between two separate diagrams obey the same principles as in ribbon graphs.

\begin{figure}[htbp]
\includegraphics[scale=1]{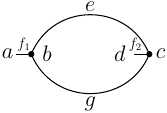}
\caption{Two-vertex ribbon graph with HRs}
\label{twovert}
\end{figure}

\begin{figure}[htbp]
   \includegraphics[scale=0.8]{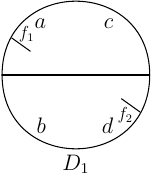}
\quad 
 \includegraphics[scale=0.8]{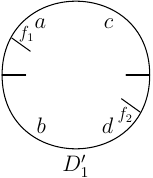}
\quad
\includegraphics[scale=0.8]{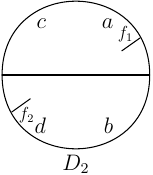}
\quad
 \includegraphics[scale=0.8]{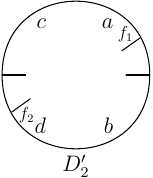}
  \caption{Related chords diagrams $D_1$, $D'_1$, $D_2$, $D'_2$}
  \label{oprotate}
\end{figure}

\begin{figure}[htbp] 
   \includegraphics[scale=0.8]{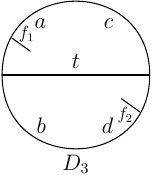}\quad
 \includegraphics[scale=0.8]{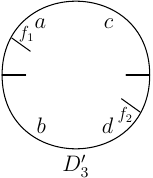}\quad
\includegraphics[scale=0.8]{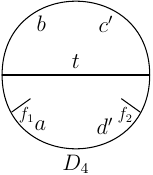}\quad
 \includegraphics[scale=0.8]{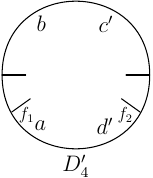}
  \caption{Related chords diagrams $D_3$, $D'_3$, $D_4$, $D'_4$}
  \label{optwiste}
\end{figure}

Consider a two-vertex ribbon graph with HRs $\cG$ with at least two edges $e$ and $g$ which are not loops. Let us write $a$, $b$, $c$ and $d$ for the sections into which $e$ and $g$ divide the cyclic orders at the vertices of $\cG$ (some HRs may be attached to the vertices as illustrated in Figure \ref{twovert}). The contractions of $e$ or of $g$ give two different bouquets with HRs. If $e$ and $g$ are positive edges, let $D_1$ be the (open) chord diagram associated with the graph we obtain by contracting $g$ in $\cG$, $D'_1$ the (open) chord diagram associated to the graph we obtain by contracting $g$ in $\cG \vee e$, $D_2$ the (open) chord diagram associated to the graph we obtain by contracting $e$ in $\cG$ and $D'_2$ the one we obtain by contracting $e$ in $\cG \vee g$ (see Figure \ref{oprotate}). If $g$ is negative (without loss of generality), we replace $D_1$, $D_2$, $D'_1$ and $D'_2$, respectively, by $D_3$, $D_4$, $D'_3$ and $D'_4$ in the previous statement (see Figure \ref{optwiste}).

In Figure \ref{optwiste}, the sector $c'$ is obtained from $c$ after a sequence of two operations:  we reverse the order of the endpoints of the HRs and chords of $c$ and we change the sign of any chord from $c$ to the rest of the diagram. The same apply to $d'$ obtained from $d$.

Two signed (open) chord diagrams are related by a rotation about the chord $e$ if they are related as $D_1$ and $D_2$ in Figure \ref{oprotate}, and   they are related by a twist about $e$, if they are related as $D_3$ and $D_4$ in Figure \ref{optwiste}. Now we can give the definitions of $R$-equivalent diagrams and the sum of two chord diagrams.  

\begin{definition}[$R$-equivalence relation \cite{bollo}]
Two diagrams or signed  diagrams $D_1$ and $D_2$ are $R$-equivalent if and only if they are related by a sequence of rotations and twists. We write  $D_1\sim D_2$.
\end{definition}
\begin{definition}[Sum of  diagrams \cite{bollo}]
The sum of two  diagrams or signed  diagrams $D_1$ and $D_2$ is obtained by choosing a point $p_i$ (not the end-point of a chord or a HR) on the boundary of each $D_i$, joining the boundary circles at these points and then deforming the result until it is again a circle.
\end{definition}
 By choosing the $p_i$ differently, this sum can be formed in many different ways but we shall show that all of them are $R$-equivalent. 
\begin{lemma}\label{lem3}
If two diagrams $D$ and $D'$ are both sums of diagrams $D_1$ and $D_2$, then they are $R$-equivalent.
\end{lemma}
\proof
The proof is the same as in \cite{bollo2} since the rotations and twists about chords move only the points $p_1$ or $p_2$ chosen on $D_1$ or $D_2$, respectively. The only fact that one must pay attention is to respect the cyclic order of the HRs on the resulting circle. In the case where there are some HRs coming  before the chord we want to rotate about or twist about, we must rotate or twist the HR about a chord before the next step.
\qed

\vspi\noindent{\bf Canonical chord diagrams.} For $i\geq 0$,  $0\leq 2j\leq i$, $0\leq k\leq i+1$, $l\geq0$ and $0\leq m\leq 2$, let $D_{i,j,k,(s;l_1,\ldots ,l_q),m}$ be the chord diagram consisting of $i$ chords, $j$ pairs of positive chords intersecting each other, $k$ connected components of the boundary graph of this diagram, $l$ HRs ($l = s + \sum_{p=1}^q l_p$) disposed in a specific way and $m$ negative chords (or twisted chords) intersecting no other chords (hence $i-2j-m$ is the number of positive chords intersecting no other chords); if $s=0$ then $q=k$, 
and if $s>0$, then $q=k-1$. 
Note that the above inequalities 
defining the canonical diagram
are not independent: if $l=0$, then $k=0$, otherwise $l>0$ and then $0<k\le i+1$. This diagram is drawn in such a way that there is a number $l-s$ of HRs partitioned in $(l_p)_{p=1,\cdots, q}$ positive chords intersecting no other chords (we shall also call these isolated chords) and $s$ is the rest of the HRs. We put ``$\tw$'' for only twisted chords for simplicity. All these chords and HRs are arranged around the circle of the diagram (see an illustration for $D_{4,1,2,(3;1),1}$ and  $D_{5,1,2,(0;1,2),1}$):
\begin{figure}[htbp]
 \includegraphics[scale=0.8]{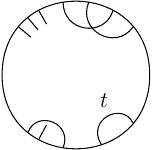}
\hspace{1cm}
\includegraphics[scale=0.8]{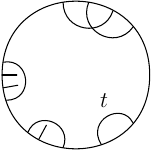}
\caption{Canonical diagrams: $D_{4,1,2,(3;1),1}$ and $D_{5,1,2,(0;1,2),1}$}
  \label{canoni}
\end{figure}

If there are no HRs on the graph, our canonical diagram corresponds exactly  to that of Bollob\'as and Riordan \cite{bollo}. Consider now a chord diagram $D$ with $l>0$ HRs. Forgetting about the HRs for a moment, one performs a sequence of rotations and twists about chords in the same way as \cite{bollo} and is led to a BR canonical diagram. The HRs in $D$ were disposed on open faces (open components) which are preserved under rotations and twists. Therefore, at the end, one adds the HRs on the resulting BR canonical diagram in order to obtain the result if the same sequence of rotations and twists about chords was performed on the initial signed chord diagram $D$ considered with HRs. The issue here is the disposition of the HRs in the BR canonical diagram. We will show however that, from the knowledge of $D$, either we can directly reconstruct the new canonical diagram or find a canonical diagram $R$-equivalent to it. 

\begin{lemma}\label{lem4}
Any (open) chord diagram $D$ is $R$-equivalent to some $D_{i,j,k,(s;l_1,\ldots, l_q),m}$. 
\end{lemma}
\proof
Let $D$ be a signed (open) chord diagram with $i$ chords, $k$ connected components of the boundary
graph of $D$ and $l$ HRs.    

Suppose $l=0$. In this case $k=0$ and $D$ is $R$-equivalent to some $D_{ijm}$ in sense of \cite{bollo}. We denote it as $D_{i,j,0,(0),m}$ since the set of partitions $(s;l_1,\cdots, l_q)$ is empty.

Assume now that $l>0$. If we forget the HRs for a moment and perform a sequence of rotations and twists about chords, we obtain that $D$ is $R$-equivalent to some $D_{ijm}$, a signed chord diagram consisting of $i$ chords, $j$ pairs of positive chords intersecting each other, $i-2j-m$ isolated positive chords and $m$ ($0\leq m\leq2$) negative isolated chords. One can add now the $l$ HRs to  $D_{ijm}$. Note that there is only one internal face which passes through all the pairs of positive chords intersecting each other and all negative chords. Then inserting HRs on this face just leads to only one connected component of the boundary graph. The remaining connected components of the boundary graph can be formed by putting a number of HRs in a certain number of isolated positive chords. Some cases have to be discussed.

\begin{figure}[htbp]
\includegraphics[scale=0.8]{CanonicD5.pdf}
\hspace{1cm}
   \includegraphics[scale=0.8]{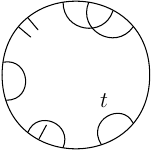}
  \caption{Two $R$-equivalent canonical diagrams: $D_{5,1,2,(0;1,2),1}\sim D_{5,1,2,(2;1),1}$}
  \label{equicanonic}
\end{figure}

Suppose at first that $i-2j-m>0$ (there is at least one positive isolated chord).
This situation decomposes in two cases. If $k\leq i-2j-m$, the number of connected components of the boundary graph is at most the number of isolated positive chords. We have two possible ways to arrange the $l$ HRs. One way is to arrange the $l$ HRs such that they are partitioned in $k$ isolated positive chords and then we obtain the canonical diagram $D_{i,j,k,(0;l_1,\ldots,l_k),m}$ ($l_p>0$, $\forall p=1,\cdots,k$). The second way is to arrange $l-s$ ($s>0$) HRs such that they are partitioned in $k-1$, $k>1$, isolated positive chords and the remaining $s$ HRs are not in any chord. Then we obtain the canonical diagram $D_{i,j,k,(s;l_1,\ldots,l_{k-1}),m}$ ($l_p>0$, $\forall p=1,\cdots,k-1$). 
If $k=1$, there is no remaining
isolated positive chords and $l=s$
yielding the canonical diagram $D_{i,j,k,(l;0),m}$. 
By a sequence of rotations  and twists about chords we have $D_{i,j,k,(0;l_1,\ldots,l_k),m}\sim D_{i,j,k,(s;l_1,\ldots,l_{k-1}),m}$ (see Figure \ref{equicanonic}.
Note
also that, for $k=1$, the above expression
trivializes to 
$D_{i,j,k,(l;0),m}\sim D_{i,j,k,(0;l),m}$.
 Now assume that $k> i-2j-m>0$, then all the $i-2j-m$ isolated positive chords of $D$ must receive some HRs. The $l-s$ ($s>0$) HRs of $D$ must be partitioned in the $i-2j-m$ chords, 
and $s$ HRs must be disposed elsewhere. 
Hence $k=i-2j-m+1$ and 
we have $D=D_{i,j,k,(s;l_1,\ldots,l_{k-1}),m}$. 

Consider finally that  $i-2j-m=0$ which means that we do not have any positive isolated chord. Then, to have a nonempty
set of HRs forces $k=1$ and then $D\sim D_{i,j,1,(l;0),m}$.
\qed

Without distinguishing 
the $s$ HRs with the remaining HRs, let us denote
$D_{i,j,k,(s; l_1,\ldots,l_{q}),m}$ as  $D_{i,j,k,(l_1,\ldots,l_{k}),m}$. Now, given a permutation $\sigma$ in $S_k$ (the permutation group with $k$ elements), 
we discussed the fact that  
$D_{i,j,k,(l_1,\ldots,l_{k}),m}\sim D_{i,j,k,(l_{\sigma(1)},\ldots,l_{\sigma(k)}),m}$. 
This simply means that the order of the sequence $(l_1,\cdots,l_k)$ does not matter
when writing the canonical diagram. 
In the following, whenever possible and for simplicity, we use $D_{i,j,k,[l],m}$ to denote $D_{i,j,k,(l_1,\cdots ,l_k),m}$.

\subsection{Universality of the polynomial $\cQ$}
We are in position to show that the 
multivariate polynomial invariant $\cQ_{\cG}$ in \eqref{brflamulti} is universal on the class of ribbon graphs with HRs. 
Then, we will prove $\mathcal{R}_{\cG}$ 
\eqref{brfla} is universal on a subclass of ribbon graphs with HRs.

We define $[l] \vdash l$ to be a partition 
of a positive integer $l \ge 0$ and, 
given $A \sset \cG$, the set of contraints on $A$ defined as: 
\bea
&&
C_{ijk[l]m}(A) \equiv 
\{ n(A)=i,\,
k(A)-F_{\inter}(A)+n(A)=j,\, C_\partial(A)=k, \,
 \, t(A)=m, f(A)=l, \cr\cr
&& 
[l] \, \text{defines the same partition as}\; [l_c] \text{ of the number of HRs of}\; A\}. 
\eea
If $C_{ijk[l]m}(A)$ holds for a given $A$, where $[l]\vdash l$
is a given partition of $l$, with parts $l_u$ ($\sum_u l_u = l$), then each $l_u$ must be a number of HRs disposed in a given connected component $\mathbf{c}_u\in \cC_\partial(A)$. 
Hence, $C_\partial(A)= k$ must be the number of parts of the partition $[l]$.

Consider the following expansion of $\cQ_{\cG}$
\bea
&&
\cQ_{\cG}(X,Y,Z,W, \{\zeta_i\})
=\sum_{i,j,k,l,m}\sum_{[l]\vdash l}\cR_{ijk[l]m}(\cG; X)   (Y-1)^{i} Z^{j} \, W^m\, \prod_{u=1}^{k}\zeta_{l_u}\, \label{brfla2},\\
&&
 \cR_{ijk[l]m}(\cG; X) :=\sum_{A\sset \cG| \;  C_{ijk[l]m}(A) \; \text{holds}} (X-1)^{\rk(\cG)-\rk(A)},
\nonumber
\eea 
where each $\cR_{ijk[l]m}(\cdot ;X)$ is  a map from the set $\mathfrak{G}$ of isomorphism classes of connected ribbon graphs with HRs to $\mathbb{Z}[X]$. 
%By equating the coefficients of $(Y-1)^{i} Z^{j} \, S^{k} \, \, W^m\prod_{u=1}^{k}\zeta_u^{l_u}$ or performing a straightforward computation from the definition of $\cR_{ijk[l]m}$, we can see from Theorem \ref{theo:BRext} that $\cR_{ijk[l]m}$ satisfies \eqref{retcondel} and \eqref{retbri}.
%For each $A\sset \cG$, a unique partition $[l] \vdash l$ contributes in \eqref{sumpart}.
%We write, 
%\bea
%&&
%\cR_{ijklm}(\cG; X)
%= \sum_{[l]\vdash l}
%\cR_{ijk[l]m}(\cG; X), 
%\label{sumpart}\\
%&&
%\cR_{ijk[l]m}(\cG; X) = 
%\sum_{A\sset \cG| 
%C_{ijklm}(A) \text{ holds and } \; l_i \text{ is the number of HRs in } \cC_i \in \cC_\partial(A)  } (X-1)^{\rk(\cG)-\rk(A)},
%\nonumber
%\eea
%We also note that
$\cR_{ijk[l]m}(\cG; X)$ fulfills the contraction-cut rules \eqref{retcondel} and \eqref{retbri} 
as given by Theorem
\ref{theo:BRext}
(the extra contraints on the type of c-subgraphs do not have any influence on the proof). 
 
Given a ring $\mathfrak{R}$ and an element $x$ of $\mathfrak{R}$, for $i$,$j$,$k$,$l$,$m$, as $\cR_{ijk[l]m}(\cdot;X)$ takes values in $\mathbb{Z}[X]$, we compose it with the ring homomorphism from $\mathbb{Z}[X]$ to $\mathfrak{R}$ mapping $X$ to $x$, and obtain a map $\cR_{ijk[l]m}(\cdot ;x)$   from $\mathfrak{G}$ to $\mathfrak{R}$. The infinite sum of these functions is of significance, but in general a finite number are non-vanishing on any given ribbon graph with HRs.

\begin{theorem}[Universality of $\cQ$]
\label{theo:univ2}
Let $\mathfrak{R}$ be a commutative ring and $x \in \mathfrak{R}$. If a function $\phi: \mathfrak{G} \to \mathfrak{R}$
satisfies 
\bea
\phi(\cG) = \left\{\begin{array}{ll} 
\phi(\cG\vee e) + \phi(\cG/e) & {\text{if e is regular}},\\\\
(x-1)\phi(\cG\vee e) + \phi(\cG/e) & {\text{if e is a bridge}},
\end{array} \right.
 \label{systemuniv2} 
\eea
\medskip 
then there are coefficients $\lambda_{ijk[l]m}\in\mathfrak{R}$, with $i\geq 0$, $0\leq k\leq i+1$, $l\geq 0$, $0\leq m\leq 1$ and $0\leq j \leq i+1$  such that 
\beq \label{univsum2}
\phi(\cG) = \sum_{i,j,k,l,m}
\sum_{[l]\vdash l}\lambda_{ijk[l]m} \cR_{ijk[l]m}(\cG;x).
\eeq 
\end{theorem}
\proof
Let us consider a two-vertex ribbon graph $\cG$ of the form in Figure \ref{twovert}. Applying equation \eqref{systemuniv2} provides two different expressions for $\phi(\cG)$: at first, one applies these relations to the positive edge $e$ and then to the positive edge $g$ (if it is not a self-loop), and then vice-versa. Equating these expressions shows that
\bea
\label{nphD1}
\phi(D_1) - \phi(D'_1)=\phi(D_2) - \phi(D'_2),
\eea
where $D_1$, $D'_1$, $D_2$ and $D'_2$ are signed chord diagrams related as illustrated in Figure \ref{oprotate}. 

Similarly, considering the case where $g$ is negative allows us to get
\bea
\label{nphD2}
\phi(D_3) - \phi(D'_3)=\phi(D_4) - \phi(D'_4),
\eea
where $D_3$, $D'_3$, $D_4$ and $D'_4$ are signed chord diagrams related as illustrated in Figure \ref{optwiste}. 

Suppose that $\phi$ satisfies \eqref{systemuniv2}, we now show that it has the form \eqref{univsum2}. We define the
$\lambda_{ijk[l]m}$  by induction. If $i=0$, then $m=0$ and we set $\lambda_{000[0]0}$ for the value of $\phi$ on a bouquet without loops and HRs, $\lambda_{011[l]0}$ ($l >0$ and $[l]= (s=l;0)$) for the value of $\phi$ on a bouquet without loops but with $l$ HRs and $\lambda_{0jk[l]m}=0$ for all other values of $j,k,l,m$.

Assume that $n\geq1$ and $\phi(\cG) =\sum_{i<n; \; j,k,l,m}
\sum_{[l]\vdash l} \lambda_{ijk[l]m} \cR_{ijk[l]m}(\cG;x)$ for all bouquets with HRs with fewer than $n$ loops. Let us set 
$\phi' = \phi -
\sum_{i<n; \; j,k,l,m}\sum_{[l]\vdash l} \lambda_{ijk[l]m} \cR_{ijk[l]m}(\cdot;x)$. $\phi'$ vanishes on bouquets with HRs with less than $n$ loops  and satisfies \eqref{systemuniv2} since $\phi$ and  the $\cR_{ijk[l]m}$ satisfy it.
This also shows $\phi'$ and $\cR_{ijk[l]m}$ obey \eqref{nphD1} and \eqref{nphD2}. 
Since $\phi'$ vanishes on chords diagrams with fewer than $n$ chords, then $\phi'(D_1) = \phi'(D_2 )$ or  $\phi'(D_3) = \phi'(D_4 )$ for related diagrams with $n$ chords. Consequently, $\phi'(D)$ depends only on the $R$-equivalence class of the chord diagram. For $j$, $k$,  $m$, $l$ and $[l]\vdash l$, there is an $\cR_{nj'k'[l']m'}$ such that $\cR_{nj'k'[l']m'}(D_{n,j'',k'',[l''],m''};x)=1$, 
if $i=i''$, $j=j''$, $k=k''$, $l=l''$, 
and $[l] = [l'']$ 
(the same partition of $l$) and 0 otherwise 
(see \eqref{phiprimlam}, \eqref{phiprimlam0} and \eqref{phiprimlam1}, in the discussion below). 
We can therefore select the $\lambda_{njk[l]m}$
 so that \eqref{univsum2} holds on the $D_{n,j,k,[l],m}$
 and this extends to 
all chord diagrams with $n$ chords.

By induction on $n$, there exist $\lambda_{ijk[l]m}$ such that \eqref{univsum2} holds for all bouquets with HRs $\cG$. The same result follows for all connected ribbon graphs with HRs using \eqref{systemuniv2}.
\qed

Let $\gamma$ be the function defined on the set $\{0,1,2\}$ by: 
\bea 
\left\{\begin{array}{ll} 
\gamma(0)=0, \\
\gamma(1)=\gamma(2)=1.
\end{array} \right. 
\eea
The computation of $\phi'$ on a canonical signed chord diagram $D_{n,j',k',[l'],m'}$ gives:
\bea
\phi'(D_{n,j',k',[l'],m'}) &=&\sum_{j,k,l,m} \sum_{[l]\vdash l} \lambda_{njk[l]m} \cR_{njk[l]m}(D_{n,j',k',[l'],m'};x)\cr &=& \sum_{j,k,l,m} \sum_{[l]\vdash l} \lambda_{njk[l]m}\delta_{j,2j'+k'+m'}\delta_{k,k'}\delta_{[l],[l']}
\delta_{m,\gamma(m')}\cr&=& \lambda_{n(2j'+k'+m')k'[l']\gamma(m')}
\label{phiprimlam}
\eea
where $\delta_{pp'}$ is $1$
if $p=p'$ for two integers $p$
and $p'$, and $0$ otherwise; the delta
function $\delta_{[l],[l']}$ of two partitions $[l]$
and $[l']$ equals 
1 if  $[l]=[l']$
(as defining the same partition of $l=l'$),
and 0 otherwise.
For some $j$, $k$, $l$ and $m$, we can compute explicitly, $\lambda_{njk[l]m}$:

$\bullet$ If $m=0$
\bea
\lambda_{njk[l]0}=\phi'(D_{n,\frac{1}{2}(j-k),k,[l],0}).
\label{phiprimlam0}
\eea 
Then $\lambda_{njk[l]0}$ is the value of $\phi'$ on the canonical signed chord diagram $D_{n,\frac{1}{2}(j-k),k,[l],0}$ if and only if $j-k\in 2\mathbb{N}$ and $j\leq n+1$. Otherwise, $\lambda_{njk[l]0}=0$.

$\bullet$ If $m=1$
\bea
\lambda_{njk[l]1}= \left\{\begin{array}{ll} 
\phi'(D_{n,\frac{1}{2}(j-k-1),k,[l],1}) & {\text{if $j-k\in 2\mathbb{N}+1$}},\\\\
\phi'(D_{n,\frac{1}{2}(j-k-2),k,[l],2}) & {\text{if $j-k\in 2\mathbb{N}+2$}}.
\end{array} \right. 
\label{phiprimlam1}
\eea

Then $\lambda_{njk[l]1}$ is the value of $\phi'$ on the canonical signed chord diagram $D_{n,\frac{1}{2}(j-k-1),k,[l],1}$ if and only if $j-k\in 2\mathbb{N}+1$ and $j\leq n+1$. It can be also  the value of $\phi'$ on the canonical signed chord diagram $D_{n,\frac{1}{2}(j-k-2),k,[l],2}$ if and only if $j-k\in 2\mathbb{N}+2$ and $j\leq n+1$. Otherwise, 
$\lambda_{njk[l]1}=0$.
	
As in case of Tutte polynomial and BR polynomial, the condition \eqref{systemuniv2} in Theorem \ref{theo:univ2} can be replaced by 
\bea \label{system2}
\phi(\cG) = \left\{\begin{array}{ll} 
\tau\phi(\cG\vee e) + \sigma\phi(\cG/e) & {\text{if e is regular}},\\\\
(x-1)\phi(\cG\vee e) + \sigma\phi(\cG/e) & {\text{if e is a bridge}},
\end{array} \right. 
\eea
with fixed element $x$, $\sigma$ and $\tau$ of $\mathfrak{R}$. If $\sigma$ and $\tau$ are invertible and $\phi(\cG)$ satisfies \eqref{system2}, then $\Phi'(\cG)=\sigma^{-\rk(\cG)} \tau^{-n(\cG)}\phi(\cG)$ satisfies \eqref{systemuniv2} with $(x-1)$ replaced by $(x-1)\sigma^{-1}$ if we want to apply Theorem \ref{theo:univ2} to this function. 

The  polynomial basis $\cR_{ijk[l]m}(\cG;-)$
cannot be easily determined from the unique knowledge of $\cR(\cG;-)$. 
We have the following expansion 
\bea
&&
\cR_{\cG}(X,Y,Z,S,W,T)
=\sum_{i,j,k,l,m}\cR_{ijklm}(\cG; X)  (Y-1)^{i} Z^{j} \, S^{k} \, W^m\, T^{l}
\crcr
&&
\cR_{ijklm}(\cG; X) 
= \sum_{[l]\vdash l}
\cR_{ijk[l]m}(\cG; X) 
\eea
where $\cR_{ijklm}(\cG; X) $ 
is a sum of the $\cR_{ijk[l]m}(\cG; X)$. 
This could be 
regarded as an obstacle 
to call $\cR(\cG;-)$ universal
on ribbon graphs with HRs. 

\begin{figure}[htbp]
\includegraphics[scale=0.6]{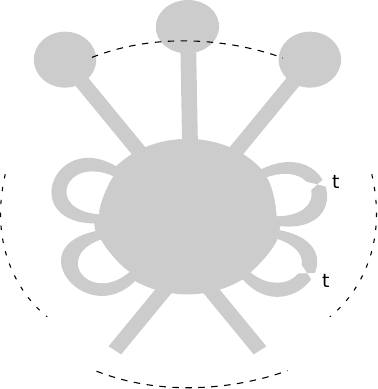}
\hspace{1cm}
\includegraphics[scale=0.6]{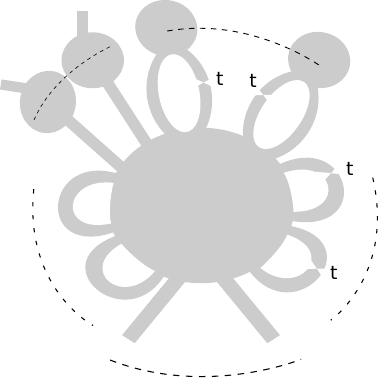}
  \caption{Ribbon graphs with HRs
  with dashed arcs corresponding to an arbitrary number of leaves, of 2-edge vertices, of (twisted) self-loops or of HRs.}
  \label{partition}
\end{figure}

A way to circumvent this issue consists in a specification of a subclass of ribbon graphs with HRs that will make  $\cR(\cG;-)$
fully characterizing a subset of  
the polynomials $\cR_{ijk[l]m}(\cG;-)$, for some precise partitions $[l]$. 
This happens when, e.g.,  $\cR_{ijklm}(\cG; X)
= \cR_{ijk[l]m}(\cG; X)$ for a specific type of partition 
$[l]$. 
One of the simplest instance where $\cR(\cG;-)$ determines $\cR_{ijk[l]m}(\cG;-)$ occurs when $k$
and $l$ fully fix the partition $[l]$
for all c-subgraphs. 
Consider a ribbon graph with the
following property: 
every connected component of the boundary of its c-subgraph contains a single HR except for one that contains the remaining $l-(k-1)$ 
HRs. Hence,  for all c-subgraphs, the partition $[l]$ is of the form $(l-(k-1), 1^{k-1})$.  This class is non empty. 
Indeed, we illustrate an example 
of one of its elements on the left of figure \ref{partition}.
Of course, the reasoning extends by requesting a fixed number, say $\alpha$, of HRs per connected component of the boundary graph of each c-subgraph, except for one having the remaining of HRs.  This is illustrated for $\alpha=2$ at the right of Figure \ref{partition}. Therein, 
one fixes $\alpha=2$ HRs per connected components of the boundary of each c-subgraphs except for one. 
Generally, 
the partition takes the form 
$[l]= (l-\alpha(k-1), \alpha^{k-1})$.
The recipe is to create a single connected component of the boundary graph of $\cG$ containing all HRs in such a way that cutting the edges keeps them in the same boundary component if they do not fall into a new connected component subgraph. 
%Indeed, consider the class of 
%single-vertex ribbon graphs will all c-subgraphs having a single connected component of their boundary graphs. Here we fix $k=1$. This class includes 
% 1-vertex ribbon graphs  with separated twisted self-loops. All c-subgraphs 
% are obtained by cutting twisted self-loops but this does not add any connected component of their boundary
% graph. 

For the subclass with $\alpha=2$, 
 we have the following reduced expansion 
\bea
\cR_{\cG}(X,Y,Z,S,W, T)
=\sum_{i,j,k,l,m}\cR_{ijk[l-2(k-1), 2^{k-1}]m}(\cG; X)  (Y-1)^{i} Z^{j} \, S^{k} \, W^m\, T^{l}.
\eea
The analog expression for any  
$\alpha \ge 1$ subclass can be easily deduced. 

The following statement becomes straightforward:  

\begin{corollary}[Universality of $\cR$]
\label{theo:univ4}
$\cR$ is universal in the sense of Theorem \ref{theo:univ2} for functions fulfilling 
the contraction/cut rule \eqref{systemuniv2}
on the class of ribbon graphs with HRs with all c-subgraphs having 2 HRs on each connected component of their boundary graph except for one which contains all remaining HRs. 
\end{corollary}
\proof 
The equation \eqref{univsum2} is now understood in the sense that
the sum over partitions $[l]$  restricts to those  
corresponding to $[l]= [l-2(k-1),2^{k-1}]$. 
The set of canonical diagrams
that appear in the present situation are restricted to those allowing precisely the considered partitions. The rest of the proof is identical to that of Theorem \ref{theo:univ2}. 
\qed 

We reasonably conjecture that the universality of $\cR$ extends beyond the $(\alpha=2)$--subclass of ribbon graphs with HRs to 
the generic $\alpha$--subclass. However, we postpone a thorough answer to the question
of the universality domain of $\cR$ to future investigation.

It is a noteworthy fact that the polynomial 
$\cR_{\cG}$ on ribbon graphs with HRs has another special  invariance. Indeed, consider 
two ribbon graphs with HRs, 
$\cG_1$ and $\cG_2$ that only differ by the
way that their HRs are distributed along the boundary components of their boundary graphs 
(in other words, the underlying ribbon graphs of $\cG_1$ and $\cG_2$ are isomorphic  and one is obtained from the other by moving around the HRs in a way of 
preserving each connected component of the
boundary graphs and the set of internal faces). In the next section, we shall give a clean definition of such operations, but for the moment, one realizes that $\cR_{\cG_1}= \cR_{\cG_2}$, for two such graphs. 
Moving around HRs on a given ribbon graph with HRs shows that there are some partitions of the set HRs in the connected components of the boundary graph that do not truly matter in the evaluation of $\cR$. 
 This strongly suggests that there exists another way of classifying ribbon graphs with HRs, a corresponding polynomial invariant that is constant on these new classes and that turns out to be universal for all maps that fulfills the same kind of invariance.

\section{Polynomial invariant for HR-equivalent  ribbon graphs}  
\label{sect:moves}
    In order to define the new category of graphs of interest, we must introduce a new equivalence relation on ribbon graphs.

\begin{figure}[htbp]
\includegraphics[scale=0.5]{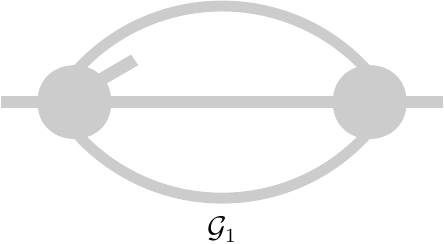}
\hspace{1cm}
   \includegraphics[scale=0.5]{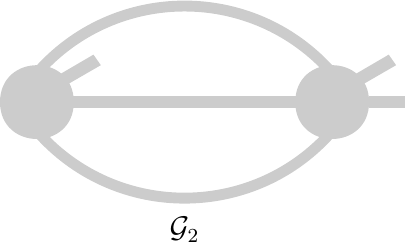}
\hspace{1cm}
 \includegraphics[scale=0.5]{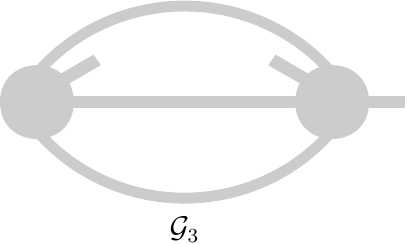}
  \caption{Some HR moves}
  \label{HR}
\end{figure}

\begin{definition}[HR move operation]
Let $\cG(\cV,\cE,\mf^0)$ be a ribbon graph with HRs. A HR move in $\cG$ consists in removing a HR $f\in\mf^0$ from one-vertex $V$ and placing $f$ either on $V$ or on another vertex such that it is called 

$\bullet$ a HR displacement if the boundary connected component
where $f$ belongs is not modified (see $\cG_1$ and $\cG_2$ in Figure \ref{HR});

$\bullet$ a HR jump if the HR is moved from one boundary connected
component to another one, provided the former remains
a connected boundary component (see $\cG_1$ and $\cG_3$ or $\cG_2$ and $\cG_3$ in Figure \ref{HR}). 
\end{definition}

One observes that under HR displacements 
the boundary graph remains unchanged whereas 
under HR jumps this graph can be modified. 
In general, under HR moves, the number of
connected components of the boundary graph 
is not modified. For instance, in Figure \ref{HR}, the graphs $\cG_2$ and $\cG_3$ are obtained from $\cG_1$
 by a HR displacement and a HR jump, respectively.

\begin{definition}[HR-equivalence relation]
We say that two ribbon graphs with HRs $\cG$ and $\cG'$ are HR-equivalent if they are related by a sequence of HR moves. 
This relation is denoted by $\cG\sim_{HR}\cG'$. 
\end{definition}
One can check that the HR-equivalence is an equivalence relation. As a consequence of the definition, 
if $\cG\sim_{HR}\cG'$, then $\cV (\cG) = \cV (\cG' ),\,$$\cE(\cG) = \cE(\cG'),\,$$\cF_{\inter}(\cG)=\cF_{\inter}(\cG')$, 
$f (\cG) = f (\cG')$, $k(\cG)=k(\cG')$, $t(\cG)=t(\cG')$, $\rk(\cG)=\rk(\cG')$, $n(\cG)=n(\cG)$ and $C_\partial{(\cG)} = C_\partial{(\cG')}$. Thus the HR moves only modify the incidence relation between HRs and vertices. We denote the HR-equivalence class of $\cG$ by $[\cG]$. Hence the three graphs in Figure \ref{HR} are HR-equivalent. For short, we will also use ``$\cG$ is equivalent to $\cG'$'' if there is no confusion. 

Let $[\cG]$ be a class of a ribbon graph with HRs under such relation. We define $V([\cG])=V(\cG),\,$$E([\cG]) = E(\cG),\,$ and $f([\cG])= f(\cG)$. The number of connected components, the rank, nullity, the number of internal faces and the number of boundary components of $[\cG]$ are those of $\cG$, namely, $k([\cG]) = k(\cG)$, $\rk([\cG]) = \rk(\cG)$, $n([\cG]) = n(\cG)$, $F_{\inter}([\cG])=F_{\inter}(\cG)$ and $C_{\partial}([\cG])=C_{\partial}(\cG)$.

The following statement holds. 

\begin{lemma}
\label{equalmon}
If two ribbon graphs with HRs $\cG$ and $\cG'$ are HR-equivalent, then  for any edge $e$ in $\cG$ and $\cG'$, $\cG\vee e$ and $\cG'\vee e$ are HR-equivalent.
\end{lemma}
\proof We shall establish that a single HR move operation commute
with cutting an edge $e$ in $\cG$. In order to do so, we must 
observe that there exists a number of connected components of the boundary graph which may 
pass through the edge $e$ and pay attention on how these components
get modified under the two processes. 

We call $\cG'$ the graph obtained from $\cG$ after the HR move. 
A case by case study is required. 

\begin{enumerate}
\item[(i)] Assume that no connected component of 
the boundary graph passes through $e$. This means that 
there is one closed face or there are two closed faces
passing through $e$. Consider in $\cG$ a HR move giving 
$\cG'$.  The HR cannot visit the closed face(s) passing through
$e$. Then after cutting  $e$, the HR cannot be hooked on the (1 or 2) boundary connected components which are generated in $\cG' \vee e$. If 
we start by cutting $e$ in $\cG$ and perform the same HR 
move in $\cG \vee e$, the HR cannot still visit the boundary components 
generated by the cut.     

\item[(ii)] Assuming now that, through $e$ pass one closed
face and one boundary connected component. Cutting
$e$ merges the close face to the boundary component. 
The reasoning is similar to the above point (i) (in the sense that 
the HR cannot be hooked to the sector generated by 
the closed face) and the operations commute. 

\item[(iii)] Let us consider now that there is no closed
face passing through $e$. Two situations, A and B, might occur: 

\medskip

A) We have a unique boundary component $C$ passing through $e$.
This case further divides into two possibilities: 

A1) The cut of $e$ generates a unique connected
component of the boundary: One easily checks that the HR
move commute with the cut. 

A2) The cut of $e$ generates two connected
components $C_1$ and $C_2$ of the boundary graph
containing each a HR coming from $e$. 

$\bullet$ Now let us assume that the move is a jump and that the HR come from another boundary component $C_0$ and ends on $C$.
After cutting $e$ that HR must be hooked to a unique $C_i$, $i=1,2$. Assuming that we cut $e$ first, the same HR jump can be
performed if and only if $C_i$ has a HR. This is  indeed
the case.

$\bullet$ Let then assume that the move is a displacement. Two situations
can happen. Either the move is done within a sector $C_i$ or 
done  from $C_1$ to $C_2$ (without loss of generality). Then we can cut $e$. If the displacement was within $C_i$, one notes that, 
after cutting $e$, we can perform the same move within the same $C_i$ which 
yields an identical configuration as above. Meanwhile, if the displacement was from $C_1$ to $C_2$ (as sectors of $C$), 
after cutting $e$, $C_1$ disconnects from $C_2$ and the same move cannot be a displacement anymore. It can be
however a jump if and only if  $C_1$ has at least one HR and this is true.  

\medskip

B) We have exactly two boundary components $C_1$
and $C_2$ passing through $e$. Note that 
the cut of $e$ generates a unique connected component $C$
of the boundary. This case divides in two further possibilities: 

$\bullet$ The move is a displacement within a sector $C_i$: 
there is no difficulty to see that the operations
commute in this case. 

$\bullet$ The move is a jump. Two further cases must be
discussed. Either the jump is from another boundary component
$C_0$ to $C_i$, $i=1,2$, then this case is again easily solved
or the jump occurs from the component $C_1$ to 
the component $C_2$ (without loss of generality).
Then, if we cut first $e$, and perform the same move,
one realizes that this move is simply a displacement
within $C$. 

\end{enumerate}

So far, we checked the case where the jump operation was defined by adding a HR to  the boundary connected 
components passing through $e$. The proof for the converse case when these components lose 
a HR after a HR jump can be done in the totally symmetric way. 
\qed

Let $[\cG]\vee e$ be the set obtained by cutting $e$ in all elements of $[\cG]$, $[\cG] - e$ the set obtained by deleting 
$e$ in all elements of $[\cG]$ and $[\cG]/e$ the set obtained by contracting $e$ in all elements of $[\cG]$. We have:

$\bullet$ $[\cG\vee e] \supset [\cG]\vee e$ and $[\cG - e] \supset [\cG] - e$.

$\bullet$ If $e$ is not a self-loop, $[\cG/e] = [\cG]/e$. 

It might happen that $[\cG]\vee e \subsetneq [\cG\vee e] $ and $[\cG] - e\subsetneq [\cG - e]$. 
Thus it is not clear that $[\cG]\vee e$ and $[\cG] - e$ correspond to some equivalence classes of some graphs.

\begin{lemma}
\label{equalpol}
 For two HR-equivalent ribbon graphs, $\cG$ and $\cG'$,  $\cR(\cG)=\cR(\cG')$ with $\cR$ the polynomial defined in \eqref{brfla}.
\end{lemma}
\proof
The proof of this lemma uses Lemma \ref{equalmon}. The number of monomials in the expansion of $\cR(\cG)$ or $\cR(\cG')$ is the same since $\cG$ and $\cG'$ have exactly the same set of edges. Each monomial of $\cR(\cG)$ is obtained from the contribution of a spanning subgraph $A\sset \cG$. Since $A$ is obtained by cutting a subset $\cE'$ of edges in $\cG$, we choose also the spanning subgraph $A'\sset \cG'$ obtained by cutting the same subset of edges in $\cG'$. Applying successively Lemma \ref{equalmon} to all elements of $\cE'$, the subgraphs $A$ and $A'$ are HR-equivalent. Then, the monomial associated with $A$ in $\cR(\cG)$ is equal to the one associated with $A'$ in $\cR(\cG')$. This achieves the proof. 
\qed

We are now ready to define the polynomial $\cR$ on HR-equivalence classes. 
\begin{definition}[Polynomial for HR-equivalence classes]
Let $\cG(\cV,\cE,\mf^0)$ be a ribbon graph with HRs  and $[\cG]$ be its HR-equivalence class. We define the polynomial of $[\cG]$ to be
\bea \label{brequifla}
\cR_{[\cG]}=\cR_{\cG}\,.
\eea
\end{definition}

The following statement is trivial.
\begin{proposition}
\label{equalpol2}
Let $\cG$ be a ribbon graph with HRs, $[\cG]$ its HR-equivalence class and $e$ one of its edges. The following relations hold 
$\cR_{[\cG\vee e]}=\cR_{\cG\vee e}$ and 
$\cR_{[\cG/e]}=\cR_{\cG/e}$.
\end{proposition}

\begin{corollary}[Contraction and cut on BR polynomial]
Let $\cG(\cV,\cE,\mf^0)$ be a ribbon graph with HRs and $[\cG]$ be its HR-equivalence class.  Then, for a regular edge $e$, 
\beq
\cR_{[\cG]}=\cR_{ [\cG\vee e]} +\cR_{[\cG/e]}\,,
\label{equiretcondel}
\eeq 
for a bridge $e$, we have 
\beq
\cR_{[\cG]} =(X-1)\cR_{[\cG \vee e]}+ \cR_{[\cG/e]}\,,
\label{equiretbri}
\eeq 
 for a trivial twisted self-loop $e$, the following holds 
\beq
\cR_{[\cG]}=\cR_{[\cG \vee e]} + (Y-1)ZW \,\cR_{[\cG/e]}\,,
\label{equiretsel2}
\eeq 
whereas for a trivial untwisted self-loop $e$, we have
\beq
\cR_{[\cG]}=\cR_{[\cG \vee e]} + (Y-1) \cR_{[\cG/e]}\,.
\label{equiretsel1}
\eeq 
\end{corollary}

\proof 
The proof of this theorem is immediate using Theorem \ref{theo:BRext} and Proposition \ref{equalpol2}.
\qed

The polynomial \eqref{brequifla} is also universal and the proof of this claim can be achieved in the same way as done for Theorem \ref{theo:univ2}. Consider the following expression:
\beq 
 \cR_{ijklm}([\cG])  := \cR_{ijklm}(\cG)\label{brfla3}
\eeq
where $\cR_{ijklm}$ keeps the meaning it has in \eqref{brfla2}. 

Consider $\mathfrak{G}_{HR}$ the set of HR-equivalence classes of isomorphism classes of connected ribbon graphs with HRs. This means that we have $\mathfrak{G}_{HR}= (\mathfrak{G}/\sim_{HR})$. Classes of chord diagrams under HR-equivalence relation are naturally 
well defined. Then the following statement holds.

Dealing with a class $[\cG]$ of $\mathfrak{G}_{HR}$, 
the important information to keep track in each partition $[l] \vdash l$ is its the number of parts that must coincide with the number of 
connected components of the boundary graph of $[\cG]$. 

\begin{theorem}[Universality of $\cR$ on classes] 
\label{theo:univ}
Let $\mathfrak{R}$ be a commutative ring and $x \in \mathfrak{R}$. If a function $\phi: \mathfrak{G}_{HR} \to \mathfrak{R}$
satisfies 
\bea \label{system}
\phi([\cG]) = \left\{\begin{array}{ll} 
\phi([\cG\vee e]) + \phi([\cG/e]) & {\text{if e is regular}},\\\\
(x-1)\phi([\cG\vee e]) + \phi([\cG/e]) & {\text{if e is a bridge}}.
\end{array} \right. 
\eea
\medskip 
Then there are coefficients $\lambda_{ijklm}\in\mathfrak{R}$, with $i\geq 0$, $0\leq k\leq i+1$, $l\geq 0$, $0\leq m\leq 1$ and $0\leq j \leq i+1$  such that 
\beq \label{univsum}
\phi([\cG]) = \sum_{i,j,k,l,m} \lambda_{ijklm} \cR_{ijklm}(x).
\eeq
\end{theorem}

\proof  
As the partition $[l]$ of HRs needs not record in HR-equivalence classes, 
we simply define the canonical diagram $\mathcal{D}_{i,j,k,l,m}$ to be $[D_{i,j,k,[l],m}]$, namely the HR-equivalence class of the canonical diagram $D_{i,j,k,[l],m}$. 
Indeed, as
$D_{i,j,k,(s; l_1, l_2, \dots, l_q),m}\sim_{HR} D_{i,j,k,(l-q; 1, 1,  \dots, 1),m}$, where
$k=q$, if $s=0$, 
or $k=q+1$ otherwise,
we do not need to track the partition of HRs in isolated positive chords in the class $\mathcal{D}_{i,j,k,l,m}$.

We adjust  the proof of Theorem \ref{theo:univ2}: if $D_1$, $D'_1$, $D_2$ and $D'_2$ are signed chord diagrams related as in Figures \ref{oprotate} or \ref{optwiste},
\bea
\phi(D_1)-\phi(D'_1)=\phi(D_2)-\phi(D'_2),
\eea
 then we can write 
\bea
\phi([D_1])-\phi([D'_1])=\phi([D_2])-\phi([D'_2]).
\eea
The definition of $\cR_{ijklm}(\cG; X)$ \eqref{brfla2} remains valid. Its expansion in terms of a sum over
partitions becomes irrelevant  
since all c-subgraphs obeying the constraints $C_{ijklm}$ are necessarily HR-equivalent.
Following step by step the proof of Theorem \ref{theo:univ2}, one proves the existence of the coefficients $\lambda_{ijklm}$ so that \eqref{univsum} holds on the $\mathcal{D}_{i,j,k,l,m}$ and therefore on all chord diagrams with $i$ chords. The rest of the proof is similar to what was done for Theorem \ref{theo:univ2}.
\qed

\vspi 
It is natural to find the restricted polynomial $\cR'$ \eqref{rpr} over classes of HR-equivalent ribbon graphs and to show its universality. 
Several other interesting developments can be now undertaken from the polynomial invariants treated in this paper. For instance, the polynomial $\cR$ does not satisfy the ordinary factorization property under the one-point-join operation (see Proposition \ref{opBRfla}). Therefore, finding a recipe theorem in the sense of \cite{joan} becomes a nontrivial task for ribbon graphs with HRs. This certainly deserves to be investigated. Furthermore, significant progresses around matroids  \cite{Duchamp:2013joa} and  Hopf algebra techniques \cite{Duchamp:2013pha} applied to the Tutte polynomial have been recently highlighted. These studies should find as well an extension for the present types of invariants. Finally, combining 
some ideas of this work and Hopf algebra calculations \cite{Raasakka:2013kaa}, another important investigation would be 
to find a universality theorem for polynomial invariants over stranded graphs \cite{rca} extending ribbons with HRs.

\begin{center}
{\bf Acknowledgements}
\end{center}

{\footnotesize 
The authors would like to thank the reviewers for their insightful comments and contributions leading to the improvement of our manuscript.
RCA research is partially supported by the Alexander von Humboldt foundation. This research was supported in part by Perimeter Institute for Theoretical Physics. Research at Perimeter Institute is supported by the Government of Canada through Industry Canada and by the Province of Ontario through the Ministry of Research and Innovation. This work is partially supported by the Abdus Salam International Centre for Theoretical Physics (ICTP, Trieste, Italy) through the Office of External Activities (OEA)-Prj-15. The ICMPA is also in partnership with the Daniel Iagolnitzer Foundation (DIF), France.}

\vspace{0.3cm}

\begin{center}
\rule{3cm}{0.01cm}
\end{center}

\end{document}